\UseRawInputEncoding
\documentclass[11pt]{article}
\usepackage{amsmath}
\usepackage{amsmath,amssymb,amsthm,graphicx,mathrsfs}
\usepackage[numbers,sort&compress]{natbib}
\usepackage{color}
\usepackage[colorlinks,bookmarksopen,bookmarksnumbered,citecolor=black, linkcolor=blue, urlcolor=black]{hyperref}
\numberwithin{equation}{section}
\newcommand{\beq}{\begin{equation}}
\newcommand{\enq}{\end{equation}}

\newtheorem{Theorem}{Theorem}[section]
\newtheorem{Lemma}[Theorem]{Lemma}

\newtheorem{Definition}[Theorem]{Definition}
\newtheorem{Remark}[Theorem]{Remark}

\newcommand{\benu}{\begin{enumerate}}
\newcommand{\beqa}{\begin{eqnarray}}
\newcommand{\beqan}{\begin{eqnarray*}}
\newcommand{\eay}{\end{array}}
\newcommand{\edm}{\end{displaymath}}
\newcommand{\eenu}{\end{enumerate}}
\newcommand{\eeq}{\end{equation}}
\newcommand{\eeqa}{\end{eqnarray}}
\newcommand{\eeqan}{\end{eqnarray*}}

\newcommand{\br}{\begin{Remark}}
\newcommand{\er}{\end{Remark}}

\newcommand{\bqa}{\begin{eqnarray}}
\newcommand{\eqa}{\end{eqnarray}}
\newcommand{\bqw}{\begin{eqnarray*}}
\newcommand{\eqw}{\end{eqnarray*}}

\newcommand{\non}{\nonumber}
\newcommand{\bea}{\begin{array}{cc}}
\newcommand{\ena}{\end{array}}

\allowdisplaybreaks[4]
\begin{document}
\begin{center}

{\large \bf Existence and regularity of global attractors for a Kirchhoff wave equation with strong damping and memory}\\

\vspace{0.20in}Bin Yang$^{1}$ $\ $ Yuming Qin $^{2,\ast}$ $\ $ Alain Miranville $^{3}$ $\ $ Ke Wang $^{2}$\\
\end{center}
$^{1}$ College of Information Science and Technology, Donghua University, Shanghai, 201620, P. R. China.\\
$^{2}$ Department of  Mathematics, Institute for Nonlinear Science, Donghua University, Shanghai, 201620, P. R. China.\\
$^{3}$ Laboratoire de Math\'ematiques et Applications, Universit\'e de Poitiers, UMR CNRS 7348-SP2MI, Boulevard Marie et Pierre Curie-T\'el\'eport 2, F-86962, Chasseneuil Futuroscope Cedex, France.\\
 \vspace{3mm}

\begin{abstract}
This paper is concerned with the existence and regularity of global attractor $\mathcal A$ for a Kirchhoff wave equation with strong damping and memory in the weighted time-dependent spaces $\mathcal H$ and $\mathcal H^{1}$, respectively. In order to obtain the existence of $\mathcal A$, we mainly use the energy method in the priori estimations, and then verify the asymptotic compactness of the semigroup by the method of contraction function. Finally, by decomposing the weak solutions into two parts and some elaborate calculations, we prove the regularity of $\mathcal A$.
\end{abstract}

\hspace{4mm}{\bf Keywords:} Kirchhoff wave equation; Global attractor; Strong damping; Memory; Regularity.

\hspace{4mm}{\bf 2020 MSC:} 35B40, 35B41, 35B65, 35L05.

\section{Introduction}
\setcounter{equation}{0}

\let\thefootnote\relax\footnote{*Corresponding author: yuming$\_$qin@hotmail.com}
\let\thefootnote\relax\footnote{\footnotesize E-mails: binyangdhu@163.com, Alain.Miranville@math.univ-poitiers.fr, kwang@dhu.edu.cn}
\quad
Let $\Omega$ be a bounded domain in $\mathbb R^{n}$ $(n \in \mathbb{N^{+}})$ with smooth boundary $\partial\Omega$. In this paper, we will discuss the following Kirchhoff wave equation with strong damping and memory
\begin{equation}
\left\{\begin{array}{ll}
\varepsilon(t)u_{tt} -M(\|\nabla u\|^2) \Delta u-\Delta u_{t}-k(0) \Delta u-\int_0^{\infty} k^{\prime}(s) \Delta u(t-s) d s
+\delta f(u_{t})\\=g(u)+h(x) & in\enspace \Omega \times\mathbb{R^{+}}, \\
u(x,t)=0 & on \enspace\partial \Omega\times\mathbb{R^{+}}, \\
u(x, 0)=u_{0}(x), \,\,u_{t}(x,0)=u_{1}(x) &\, x \in \Omega,
\end{array}\right.\label{1.1-5}
\end{equation}
where $\varepsilon(t)$ is a time-dependent function, $M(\cdot)$ is a nonlocal function satisfying $M(s)=s^{\frac{m}{2}}$ for any $m \geqslant 1$ and $s \in \mathbb{R}^{+}$, $-\Delta u_{t}$ is a strong damping, $k(\cdot)$ admits $k(0)>0$, $k(\infty)>0$ and $k^{\prime}(s) \leqslant 0$, $\delta f(u_{t})$ is a perturbation term with $0<\delta<1$, $g(u)$ and $h(x) \in L^{2}(\Omega)$ are both nonlinear functions.

The time-dependent function $\varepsilon(t) \in C^{1}(\mathbb{R})$ is decreasing and satisfying
\begin{equation}
\lim _{t \rightarrow+\infty} \varepsilon(t)=a,
\label{1.2-5}
\end{equation}
for any $0<a\ll1$ and there exists a constant $L>a>0$ such that
\begin{equation}
\sup _{t \in \mathbb{R}}(|\varepsilon(t)|+|\varepsilon^{\prime}(t)|) \leqslant L.
\label{1.3-5}
\end{equation}

Assume $\delta f\left(u_t\right) \in C^{1}(\mathbb{R})$ and $f(\cdot)$ satisfies
\begin{equation}
f(0)=0, \label{1.4-5}
\end{equation}
\begin{equation}
\liminf _{|s| \rightarrow \infty} f^{\prime}(s)>0
\label{1.5-5}
\end{equation}
and
\begin{equation}
|f(s)| \leqslant C_1(1+|s|^{p_1}),
\label{1.6-5}
\end{equation}
where $p_{1}=0$ and $p_{1}\in\left(0, \frac{n+2}{n-2}\right)$ for $n=1, 2$ and $n \geqslant 3$, respectively.

Furthermore, suppose $g(u) \in C^{1}(\mathbb{R})$ and there exists a constant $p_{2}$ that fulfills $p_{2}=1$ when $n=1, 2$ and $p_{2}\in\left(1, \frac{n+2}{n-2}\right)$ when $n \geqslant 3$ such that
\begin{equation}
\left|g^{\prime}(s)\right| \leqslant C_2\left(1+|s|^{p_2-1}\right)
\label{1.7-5}
\end{equation}
and
\begin{equation}
\lim _{|s| \rightarrow \infty} \inf \frac{g(s)}{s}<\lambda_1,
\label{1.8-5}
\end{equation}
where $\lambda_1>0$ is the constant in the Poincar\'{e} inequality $\lambda_1\|u\|^2 \leqslant \|\nabla u\|^2 $. In addition, let $G(u)=\int_0^u g(u) d s$.

As in \cite{gmp.5}, we define
\begin{equation}
\eta^t(x, s)=u(x, t)-u(x, t-s).
\label{1.9-5}
\end{equation}

Taking the derivative of $(\ref{1.9-5})$ with respect to $x$ and $t$, assuming $\mu(s)=-k^{\prime}(s)$ and $k(\infty)=1$, then we can transform $(\ref{1.1-5})$ as
\begin{equation}
\varepsilon(t) u_{t t}-\left(1+\|\nabla u\|^m\right) \Delta u-\Delta u_t-\int_0^{\infty} \mu(s) \Delta  \eta^t(x, s) d s+\delta f\left(u_t\right)=g(u)+h(x)
\label{1.10-5}
\end{equation}
and
\begin{equation}
\eta_t=-\eta_s+u_t
\label{1.11-5}
\end{equation}
with initial boundary value conditions
\begin{equation}
\begin{cases}u(x, t)=0, & x \in \partial \Omega, \,t \geqslant 0, \\
 \eta^t(x, s)=0, & (x, s) \in \partial \Omega \times \mathbb{R}^{+}, \,t \geqslant 0, \\
 u(x, 0)=u_0(x), & x \in \Omega, \\
 u_t(x, 0)=u_1(x), & x \in \Omega, \\
 \eta^t(x, 0)=0, & x\in\Omega,\\
 \eta^0(x, s)=\eta_0(x, s), & (x, s) \in \partial \Omega \times \mathbb{R}^{+}.
 \end{cases}
 \label{1.12-5}
\end{equation}

The function $\mu(\cdot)$ in (\ref{1.10-5}) is a memory kernel which satisfies the following assumptions
\begin{equation}
\mu (s)\in C^1\left(\mathbb{R}^{+}\right) \cap L^1\left(\mathbb{R}^{+}\right),
 \label{1.13-5}
\end{equation}
\begin{equation}
\mu(s) \geqslant 0,\, \mu^{\prime}(s) \leqslant 0,
 \label{1.14-5}
\end{equation}
\begin{equation}
\int_0^{\infty} \mu(s) \mathrm{d} s=\delta_1>0,
 \label{1.15-5}
\end{equation}
\begin{equation}
\mu^{\prime}(s)+\delta_2 \mu(s) \leqslant 0,
 \label{1.16-5}
\end{equation}
for any $s \in \mathbb{R}^{+}$, where $\delta_1$ and $\delta_2$ are independent of the parameter $\delta$ in problem $(\ref{1.1-5})$ and satisfying $\delta_2\geqslant\delta_1>0$.

Next, we will introduce the definitions and properties of some useful spaces in this paper. %and assumptions of the functions in problem $(\ref{1.1-5})$.
Suppose $H_{1}=L^{2}(\Omega)$, $H_{2}=H_{0}^{1}(\Omega)$, and $H_{3}=H^{2}(\Omega)$, and denote their norms as $\|\cdot\|_{H_{1}}=\|\cdot\|$, $\|\cdot\|_{H_{2}}=\|\nabla\cdot\|$ and $\|\cdot\|_{H_{3}}=\|\Delta\cdot\|$, respectively.

For the relative displacement history, we assume $\mathcal M_{i}$ $(i=1, 2, 3)$ are $L^{2}_{\mu}$-weighted Hilbert spaces, which satisfy
\begin{equation}
\mathcal M_i=L_\mu^2\left(\mathbb{R}^{+} ; H_i\right)=\left\{\xi(s):\mathbb{R}^{+} \rightarrow H_i \, \,|\, \|\xi(s)\|_{\mathcal M_i}^2<\infty\right\}
\label{1.17-5}
\end{equation}
and their inner products and norms are defined as
\begin{equation}
(\xi(s), \eta(s))_{\mathcal M_i}=\int_0^{\infty} \mu(s)(\xi(s), \eta(s))_{H_i} d s
\label{1.18-5}
\end{equation}
and
\begin{equation}
\|\xi(s)\|_{\mathcal M_i}^2=\int_0^{\infty} \mu(s)\|\xi(s)\|_{H_i}^2 d s,
\label{1.19-5}
\end{equation}
respectively.

Moreover, our phase space is the weighted time-dependent space $\mathcal{H}=H_{2}\times H_{1}\times \mathcal{M}_{2}$, which is equipped with the norm
\begin{equation}
\|(u, u_t, \xi)\|_{\mathcal H}^2=\|\nabla u\|^2+\varepsilon(t)\|u_t\|^2+\|\xi\|_{\mathcal M_2}^2.
\label{1.20-5}
\end{equation}

Additionally, the weighted time-dependent space $\mathcal{H}^{1}=H_{3}\times H_{2}\times \mathcal{M}_{3}$ is more regular than $\mathcal H$ and is endowed with the norm
\begin{equation}
\left\|\left(u, u_t, \xi\right)\right\|_{\mathcal{H}^{1}}^2=\|\Delta u\|^2+\varepsilon(t)\|\nabla u_t \|^2+\|\xi\|_{\mathcal{M}_3}^2.
\label{1.21-5}
\end{equation}

The Kirchhoff type wave equations were widely discussed after Kirchhoff \cite{k.5} first studied the oscillation of elastic strings for a Kirchhoff type strong oscillatory wave equation $u_{tt}-\left(1+\epsilon_1\|\nabla u\|^2\right) \Delta u+f(u)=g(x)$ with $\epsilon_1>0$. Zhu \cite{z.5} got the existence of global attractor for the  Kirchhoff type wave equation $u_{t t}-M\left(\|\nabla u\|^2\right) \Delta u-\alpha \Delta u_{t t}-\Delta u_t+f(u)=h(x)$ in $H^1(\Omega) \times L^2(\Omega)$. Besides, Yang \cite{yzj.5} obtained the Kirchhoff type equation $u_{t t}-M\left(\|\nabla u\|^2\right) \Delta u-\Delta u_t+u+u_t+g(x, u)=f(x)$ has a global attractor in $H^2(\Omega) \times H^1(\Omega)$, which is connected and has finite fractal and Hausdorff dimension. Then Yang and Wang \cite{yw.5} proved the existence of the semigroup for the Kirchhoff type equation with a strong dissipation $u_{t t}-M\left(\|\nabla u\|^2\right) \Delta u-\Delta u_t+$ $h\left(u_t\right)+g(u)=f(x)$ in $V_{1+\rho} \times V_\rho$ with $0<\rho \leqslant 1$, where $V_\rho=D(A^{\frac{\rho}{2}})$ is a Hilbert space with the scalar product $(u, v)_\rho=(A^{\frac{\rho}{2}} u, A^{\frac{\rho}{2}} v)$ and the norm $\|u\|_{V_\rho}=\|A^{\frac{\rho}{2}} u\|$ and $A$ is an operator satisfying $A u=-\Delta u$ for any $u \in H^2(\Omega) \cap H_0^1(\Omega)$. Later, Yang and Li \cite{yl.5} discussed the existence of the finite-dimensional global attractors and exponential attractors for the same equation and in the same phase space as \cite{yw.5}. Furthermore, Wu \cite{w.5} verified the solution $u$ of the nonlinear Kirchhoff type wave equation $u_{t t}-M\left(\|\nabla u\|_2^2\right) \Delta u+h(t) g\left(u_t\right)+f(u)=0$ satisfies $E(t) \rightarrow 0$ as $t \rightarrow \infty$ in $H_0^1(\Omega)\times L^2(\Omega)$, where $E(t)=\frac{1}{2}\left\|u_t\right\|_2^2+\frac{1}{2} \bar{M}\left(\|\nabla u(t)\|_2^2\right)+\int_{\Omega} F(u) d x$ with $\bar {M}(\lambda)=\int_0^\lambda M(s) d s$, $F(u)=\int_0^u f(s) d s$ and $\lambda>0$.

In addition, Wang and Zhong \cite{wz.5} investigated the upper semicontinuity of pullback attractors for the Kirchhoff wave model with strong damping $u_{tt}-\Delta u_t-\left(1+\epsilon_2 \|\nabla u\|^2\right) \Delta u+$ $f(u)=g(x, t)$ with $\epsilon_2 >0$ in $H_0^1(\Omega) \times L^2(\Omega)$. Li and Wei \cite{lw.5} considered the existence of the global attractor of the quasilinear Kirchhoff type wave equation $\bar{\mu} u_{t t}-\left(1+\|\nabla u\|^2\right) \Delta u+u_t+f(u)=g(x)$ with $\bar{\mu}>0$ in $V_2$. Later on, Yang and Ding \cite{ydi.5} established the well-posedness, the existence of global and exponential attractors in critical nonlinearity case for the Kirchhoff equation $u_{t t}-\Delta u_t-M\left(\|\nabla u\|^2\right) \Delta u+u_t+g(x, u)=f(x)$ in $H^1\left(\mathbb{R}^N\right) \times L^2\left(\mathbb{R}^N\right)$ with $N\geqslant3$. Li \cite{l.5} proved the existence of the global attractor of the Kirchhoff type equation with a strong dissipation $u_{t t}-\left(1+\|\nabla u\|^2\right) \Delta u-\Delta u_t+f\left(u_t\right)+g(u)=h(x)$ in $(H^2(\Omega) \cap H_0^1(\Omega)) \times H_0^1(\Omega)$. Moreover, Li and Yang \cite{ly.5} analyzed the robustness of pullback attractors and pullback exponential attractors for the Kirchhoff wave model with strong nonlinear damping $u_{t t}-\left(1+\epsilon_3\|\nabla u\|^2\right) \Delta u-\sigma\left(\|\nabla u\|^2\right) \Delta u_t+f(u)=g(x, t)$ in $(H_0^1(\Omega) \cap L^{p+1}(\Omega))\times L^2(\Omega)$ with $1 \leqslant p \le \frac{N+4}{(N-4)^{+}}$ and $N\geqslant3$, where $\epsilon_3 \in[0,1\rceil$ is an extensibility parameter, $\sigma \in C^1\left(\mathbb{R}^{+}\right)$ and $\sigma(s) \geqslant \sigma_0>0$. Yang and Da \cite{yd.5} learned the existence of global and exponential attractors of the Kirchhoff wave equation $u_{t t}-\left(1+\epsilon_4\|\nabla u\|^2\right) \Delta u-\Delta u_t+h\left(u_t\right)+$ $g(u)=f(x)$ with $\epsilon_4 \in[0,1]$ in $H_0^1(\Omega)\times L^2(\Omega)$. Additionally, some scholars also considered the long time behavior of solutions to the Kirchhoff type equations (see \cite{c.5, ljs.5, lyf.5, mz.5}).

Next, we will point out the difficulties and innovations of this paper.

(1) As far as we know, there are about 40 published works on the attractors of Kirchhoff wave equations at present. Obviously, the results are less than those of other types of wave equations. This is due to the fact that the Kirchhoff term $M(\cdot)$ in the Kirchhoff wave equation is nonlocal, which greatly increases the difficulty and the amount of calculations and estimations in studying the long time behavior of the solutions. To simplify the calculations, we first convert problem (\ref{1.1-5}) into problem $(\ref{1.10-5})-(\ref{1.12-5})$ by using the method in \cite{gmp.5}. Then in the energy estimations, we skillfully set up hypotheses $(\ref{1.15-5})$ and $(\ref{1.16-5})$, and through lots of delicate estimations, we finally obtain the desired results.

(2) In addition, our phase space $\mathcal H$ is a weighted time-dependent space and its norm contains a time-dependent function $\varepsilon(t)$,  which makes our problem $(\ref{1.10-5})-(\ref{1.12-5})$ closer to the real physical world. After considering the existence of global attractor $\mathcal A$, we also verify its regularity in the weighted time-dependent space $\mathcal H^1$ by a series of energy estimates. For the Kirchhoff wave equations, there is \textbf{no any result} on the existence of attractors in a weighted time-dependent space until now. Our paper is a brand new attempt, which will give some new inspiration to explore the existence and regularity of attractors for the Kirchhoff wave equations in weighted time-dependent spaces.

(3) Last but not least, our problem $(\ref{1.10-5})-(\ref{1.12-5})$ includes not only the strongly damped term $-\Delta u_{t}$, but also the memory term $-\int_0^{\infty} k^{\prime}(s) \Delta u(t-s) d s$ and the disturbance term $\delta f(u_{t})$, which are more comprehensive than those papers that only consider the strongly damped term or the memory term or the disturbance term. Meanwhile, these terms ensure our problem has a wide range of applications.

This paper is organized as follows. In $\S 2$, we first introduce some basic definitions and lemmas. Then we verify in $\S 3$ the existence of global attractor $\mathcal A$ in the weighted time-dependent space $\mathcal H$ by the methods of energy estimation and contraction function. Finally, using some delicate calculations, we obtain the regularity of $\mathcal A$ in the weighted time-dependent space $\mathcal H^{1}$ in $\S 4$.
\section{Preliminaries}
$\ \ \ \ $ In this section, we will introduce some definitions and lemmas related to global attractors.

Let $X$ be a Banach space or a closed subset of a Banach space.

\begin{Definition} {\rm(\cite{r})}
If $\mathcal B_t(R)$ is centered at the origin and satisfies
$$
\mathcal B_t(R)=\left\{u \in X : \|u\|_{X} \leqslant R\right\}
$$
for any constant $R>0$, then $\mathcal B_t(R)$ is called a $R$-ball in $X$.
\label{def2.1-5}
\end{Definition}

\begin{Definition} {\rm(\cite{r})}\label{def2.2-5}
The Hausdorff semi-distance of two nonempty  sets $A, B \subset X$ is defined as
$$
dist_{X}(A, B)=\sup _{x \in A} \inf _{y \in B}\|x-y\|_{X} \, .
$$
\end{Definition}

\begin{Definition}{\rm(\cite{qs.5,r})} A parameter family $S(t): X \rightarrow X$ is called a semigroup if

(i) $S(0)=Id$ is the identity operator in $X$;

(ii) $S(t+s)=S(t)S(s)$, for any $s, t \geqslant 0$.

\label{def2.3-5}
\end{Definition}

\begin{Definition}{\rm(\cite{h.5,r})}\label{def2.4-5} A subset $\mathcal{A}$ in $X$ is called a global attractor if and only if

(i) $\mathcal{A}$ is invariant, i.e., $S(t) \mathcal{A}=\mathcal{A}$ for all $t \geq 0$;

(ii) $\mathcal{A}$ is compact in $X$;

(iii) $\mathcal{A}$ attracts every bounded subset $B$ in $X$, that is,$$\lim _{t \rightarrow +\infty} dist_{X}\left(S(t) B, \mathcal A\right)=0.$$
\end{Definition}

\begin{Definition}{\rm(\cite{h.5,r})}
A bounded set $D \subset X$ is called a absorbing set, if there exists a $T>0$ such that $S(t)B \subset D$ for any bounded subset $B \subset X$ and any $t \geqslant T$.
\label{def2.5-5}
\end{Definition}

\begin{Definition}{\rm(\cite{h.5})}\label{def2.6-5}
Let $X$ be a Banach space and $B$ be a bounded subset of $X$, we call a function $\Phi(\cdot, \cdot)$ which defined on $X \times X$, is a contractive on $B \times B$ if for any sequence $\left\{x_n\right\}_{n=1}^{\infty} \subset B$, there exists a subsequence $\left\{x_{n_k}\right\}_{k=1}^{\infty} \subset\left\{x_n\right\}_{n=1}^{\infty}$ such that
$$
\lim _{k \rightarrow \infty} \lim _{l \rightarrow \infty} \Phi_T\left(x_{n_k}, x_{n_l}\right)=0.
$$
Moreover, let $C(B)$ represent all contractive functions on $B \times B$.
\end{Definition}

\begin{Lemma}{\rm(\cite{h.5})}
Let $\{S(t)\}_{t \geq 0}$ be a semigroup in $X$ and it has a bounded absorbing set $B_0$. Furthermore, assume that for any constant $\sigma>0$, there exists $T(B_0, \sigma) \leqslant t$ and $\Phi(\cdot, \cdot) \in C(B)$ such that
$$
\|S(T) x-S(T) y\|_{X_t} \leqslant \sigma+\Phi_T(x, y), \quad \forall\, x, y \in B_0,
$$
then $\{S(t)\}_{t \geq 0}$ is asymptotically compact in $X$, i.e., for any bounded sequence $\left\{y_n\right\}_n^{\infty} \subset X$ and $\left\{t_n\right\}$ with $t_n \rightarrow +\infty,\left\{S\left(t_n\right) y_n\right\}_{n=1}^{\infty}$ is compact in $X$.
\label{lem2.7-5}
\end{Lemma}

\section{Existence of global attractor $\mathcal A$}
$\ \quad$ In this section, we will show the existence of global attractor $\mathcal A$. First of all, by using the Faedo-Galerkin approximation method, we can easily establish the following theorem about the existence, uniqueness and continuity to the weak solutions for problem $(\ref{1.10-5})-(\ref{1.12-5})$ in the weighted time-dependent space $\mathcal H$.

\begin{Theorem}
Under the assumptions of $M(\cdot)$, $\varepsilon(t)$, $k(\cdot)$, $\mu(\cdot)$, $\delta$, $f$, $g$ and $h$ in $\S 1$, if $z_{0}=(z_{0}, (u_{0})_t, \eta_{0}) \in \mathcal H_{0}$, then there exists an unique weak solution $z=\left(u, u_t, \eta\right)$ to problem $(\ref{1.10-5})-(\ref{1.12-5})$ which satisfies $\left.z \in C([0, T]; \mathcal H\right)$ for any $T>0$ in $\mathcal H$.
\label{th3.1-5}
\end{Theorem}

\begin{Lemma} Under the assumptions of $M(\cdot)$, $\varepsilon(t)$, $k(\cdot)$, $\mu(\cdot)$, $\delta$, $f$, $g$ and $h$ in $\S 1$, the semigroup $\{S(t)\}_{t \geqslant 0}$ for problem $(\ref{1.10-5})-(\ref{1.12-5})$ has a bounded absorbing set in $\mathcal H$.
\label{lem3.1-5}
\end{Lemma}
$\mathbf{Proof.}$ Multiplying $(\ref{1.10-5})$ by $u_t$ and using $G(u)=\int_0^u g(u) d s$ and $(\ref{1.11-5})$, we arrive at
%\begin{equation}
\begin{align}
&\frac{d}{d t}\left(\frac{1}{2} \varepsilon(t)\|u_t\|^2+\frac{1}{2}\|\nabla u\|^2+\frac{1}{m+2}\|\nabla u\|^{m+2}+\frac{1}{2}\|\eta \|_{\mathcal M_2}^2-\int_{\Omega}(G(u)+h u) d x\right)\nonumber\\
&-\frac{1}{2} \varepsilon^{\prime}(t) \| u_t\|^2+\| \nabla u_t \|^2+\left(\eta, \eta_s\right)_{\mathcal M_2}+\delta\left(f\left(u_t\right), u_t\right)=0.\label{3.1-5}
\end{align}
%\end{equation}

From $(\ref{1.16-5})$ and the definition of $\mathcal M_2$ , we obtain
\begin{align}
(\eta, \eta_s)_{\mathcal M_2}&=\int_0^{\infty} \mu(s) \int_{\Omega} \nabla \eta \nabla \eta_s d x d s=-\frac{1}{2} \int_0^{\infty} \mu^{\prime}(s)\|\nabla \eta\|^2 d s\non\\
& \geqslant \frac{\delta_2}{2} \int_0^{\infty} \mu(s)\|\nabla \eta\|^2 d s=\frac{\delta_2}{2}\|\eta \|^2 _{\mathcal M_2}.\label{3.2-5}
\end{align}

Let \begin{equation}
I_1(t)=\frac{1}{2} \varepsilon(t)\|u_t\|^2+\frac{1}{2}\| \nabla u\|^2+\frac{1}{m+2}\| \nabla u\|^{m+2}+\frac{1}{2}\|\eta\|^2_{\mathcal M_2}-\int_{\Omega}(G(u)+h u) d x.
\label{3.3-5}
\end{equation}

Then inserting $(\ref{3.2-5})$ and $(\ref{3.3-5})$ into $(\ref{3.1-5})$, we conclude
\begin{equation}
\begin{aligned}
\frac{d}{d t} I_1(t)-\frac{1}{2} \varepsilon^{\prime}(t)\|u_t\|^2+\|\nabla u_t\|^2+\frac{\delta_2}{2}\|\eta\|^2_{\mathcal M_2}+\delta\left(f\left(u_t\right), u_t\right)
\leqslant 0.
\end{aligned}
\label{3.4-5}
\end{equation}

Integrating $(\ref{3.4-5})$ from $0$ to $t$, we derive
\begin{equation}
I_1(t)-\frac{1}{2} \int_0^t \varepsilon^{\prime}(s)\left\|u_s\right\|^2 d s+\int_0^t\left\|\nabla u_s\right\|^2 d s+\frac{\delta_2}{2} \int_0^t\|\eta\|_{\mathcal M_2}^2 d s+\delta \int_0^t(f(u_s), u_s)ds \leqslant I_1(0).
\label{3.5-5}
\end{equation}

By $(\ref{1.5-5})$, it follows that $f$ is increasing, then from $(\ref{1.4-5})$ there holds $f(s)>0$ for any $s\in \mathbb R^+$. Moreover, since $0<\delta<1$, we conclude $\delta(f(u_t), u)>0$. Noting that $\varepsilon (t)$ is decreasing and using $(\ref{1.2-5})$, we derive $-\frac{1}{2} \varepsilon^{\prime}(t)\left\|u_t\right\|^2>0$.

Consequently, we obtain
\begin{equation}
-\frac{1}{2} \int_0^t \varepsilon^{\prime}(s)\left\|u_s\right\|^2 d s+\int_0^t\left\|\nabla u_s\right\|^2 d s+\frac{\delta_2}{2} \int_0^t\|\eta\|_{\mathcal M_2}^2 d s+\delta \int_0^t\left(f\left(u_s\right), u_s\right)ds>0.
\label{3.6-5}
\end{equation}

Then $(\ref{3.5-5})$ and $(\ref{3.6-5})$ imply that for any $t>0$,
\begin{equation}
I_1(t)\leqslant I_1(0).
\label{3.7-5}
\end{equation}

Using $(\ref{1.8-5})$, we derive that there exists constants $\lambda$ satisfies $\lambda_1<\lambda<\frac{3}{2}\lambda_1$ and $C_3>0$ such that
\begin{equation}
(g(u), u)<\lambda_1\|u\|^2<\frac{\lambda_1}{2}\|u\|^2+C_3 |\Omega|
\label{3.8-5}
\end{equation}
and
\begin{equation}
\begin{aligned}
\int_{\Omega} G(u) d x &<\frac{\lambda_1}{2}\|u\|^2<\frac{\lambda}{4}\|u\|^2+C_3|\Omega|,
\end{aligned}
\label{3.9-5}
\end{equation}
where $|\Omega|$ is the measure of the domain $\Omega$.

Furthermore, from the Young inequality, it follows that there exists a constant $\alpha>0$ such that
\begin{equation}
-\int_{\Omega} h u d x \geqslant-\alpha\|u\|^2-\frac{1}{4 \alpha}\|h\|^2.
\label{3.10-5}
\end{equation}

Similarly, there exists constants $\lambda$ admits $0<\lambda+4 \alpha<\frac{3}{2} \lambda_1$ and $C_4>\max \left\{C_3, \frac{1}{4 \alpha}\right\}$ such that
\begin{align}
I_1(0) & \geqslant I_1(t) \geqslant \frac{1}{2} \varepsilon(t)\|u_ t\|^2+\frac{1}{8}\| \nabla u\|^2+\frac{1}{m+2}\| \nabla u\|^{m+2}+\frac{1}{2}\|\eta\|_{\mathcal M_2}^2-C_4\left(|\Omega|+\|h\|^2\right)\non\\
&\geqslant-C_4(|\Omega|+\|h\|^2).\label{3.11-5}
\end{align}

Multiplying $(\ref{1.10-5})$ by $v=u_t+\alpha u$ and using $(\ref{3.2-5})$, we derive
\begin{align}
&\frac{d}{d t}\left(\frac{1}{2} \varepsilon(t)\|v\|^2-\frac{\alpha^2}{2} \varepsilon(t)\|u\|^2+\frac{1}{2}\|\nabla u\|^2+\frac{1}{m+2}\|\nabla u\|^{m+2}+\frac{\alpha}{2}\|\nabla u\|^2+\frac{1}{2}\|\eta\|_{\mathcal{M}_2}^2\right.\non\\
&\left.-\int_{\Omega}(G(u)+h u) d x\right)-\frac{1}{2} \varepsilon^{\prime}(t)\|v\|^2-\alpha \varepsilon(t)\|u_t\|^2
+\frac{\alpha^2}{2} \varepsilon^{\prime}(t)\|u\|^2+\alpha\|\nabla u\|^2+\alpha\|\nabla u\|^{m+2}\non\\
&+\left\|\nabla u_t\right\|^2+\frac{\delta_2}{2}\|\eta\|_{\mathcal{M}_2}^2+\alpha(\eta, u)_{\mathcal{M}_2}-\alpha(g(u), u)-\alpha(h, u) \leqslant 0.\label{3.12-5}
\end{align}

Let
\begin{align}
A_1(t)&=\frac{1}{2} \varepsilon(t)\|v\|^2-\frac{\alpha^2}{2} \varepsilon(t)\|u\|^2+\frac{1}{2}\|\nabla u\|^2+\frac{1}{m+2}\|\nabla u\|^{m+2}+\frac{\alpha}{2}\|\nabla u\|^2+\frac{1}{2}\left\|\eta\right\|_{\mathcal M_2}^2\non\\
&-\int_{\Omega}(G(u)+h u) d x \label{3.13-5}
\end{align}
and
\begin{align} B_1(t)&=-\frac{1}{2} \varepsilon^{\prime}(t)\|v\|^2-\alpha \varepsilon(t) \|u_t\|^2+\frac{\alpha^2}{2} \varepsilon^{\prime}(t)\| u\|^2+\alpha\| \nabla u\|^2+\alpha\| \nabla u\|^{m+2}+\| \nabla u_t\|^2\non\\
&+\frac{\delta_2}{2}\|\eta\|_{\mathcal M_2}^2+\alpha(\eta, u)_{\mathcal M_2}-\alpha(g(u), u)-\alpha(h, u).\label{3.14-5}
\end{align}

Then from $(\ref{3.12-5})-(\ref{3.14-5})$, we deduce
\begin{equation}
\frac{d}{d t} A_1(t)+B_1(t) \leqslant 0.
\label{3.15-5}
\end{equation}

From $(\ref{3.9-5})$, we conclude
\begin{equation}
-\int_{\Omega} G(u) d x>-\frac{\lambda}{4}\|u\|^2-C_3|\Omega|.
\label{3.16-5}
\end{equation}

Besides, noting that $h \in L^2(\Omega)$ and using the Young inequality, we obtain
\begin{equation}
-\int_{\Omega} h u d x \geqslant-\frac{1}{2}\|h\|^2-\frac{1}{2}\|u\|^2.
\label{3.17-5}
\end{equation}

Thanks to $(\ref{3.16-5})$, $(\ref{3.17-5})$ and the Poincar\'{e} inequality, we derive that if the constant $\alpha$ further satisfies $-\frac{\alpha^2}{2 \lambda_1} \varepsilon(t)-\frac{\lambda}{4 \lambda_1}-\frac{1}{2 \lambda_1}+\frac{1+\alpha}{2}>\frac{1}{8}$, then we arrive at
\begin{align}
&-\frac{\alpha^2}{2} \varepsilon(t)\|u\|^2+\frac{1}{2}\|\nabla u\|^2+\frac{1}{m+2}\|\nabla u\|^{m+2}+\frac{\alpha}{2}\|\nabla u\|^2+\frac{1}{2}\|\eta\|^2_{\mathcal M_2}-\int_{\Omega} (G(u)+h u) d x\non\\
& \geqslant\frac{1}{8}\|\nabla u\|^2+\frac{1}{m+2}\|\nabla u\|^{m+2}+\frac{1}{2}\|\eta\|_{\mathcal M_2}^2-C_3| \Omega|-\frac{1}{2}\|h\|^2.\label{3.18-5}
\end{align}

Therefore, there exists a constant $C_5>\max \left\{C_3, \frac{1}{2}\right\}$ such that
\begin{equation}
A_1(t)\geqslant\frac{1}{2} \varepsilon(t)\|v\|^2+\frac{1}{8}\|\nabla u\|^2+\frac{1}{m+2}\|\nabla u\|^{m+2}+\frac{1}{2}\|\eta\|_{\mathcal  M_2}^2-C_5\left(|\Omega|+\|h\|^2\right).
\label{3.19-5}
\end{equation}

Taking the two parameters in the Young inequality as $\frac{\delta_1}{4}$ and $\frac{1}{\delta_1}$, then after some simple calculations, we conclude from (\ref{1.15-5}) that
\begin{equation}
\alpha(\eta, u)_{\mathcal M_2} \geqslant-\frac{\delta_1}{4}\|\eta\|_{\mathcal M_2}^2-\alpha^2\|\nabla u\|^2.
\label{3.20-5}
\end{equation}

By $(\ref{3.9-5})$, we obtain
\begin{equation}
-(g(u), u)>-\frac{\lambda_1}{2}\|u\|^2-C_3|\Omega|.
\label{3.21-5}
\end{equation}

Carrying out the calculations similar to $(\ref{3.10-5})$, we deduce
\begin{equation}
-\alpha(h, u) \geqslant-\frac{\alpha}{2}\|h\|^2-\frac{\alpha}{2}\|u\|^2.
\label{3.22-5}
\end{equation}

Then inserting $(\ref{3.20-5})-(\ref{3.22-5})$ into $(\ref{3.14-5})$ and noting $\delta_2\ge\delta_1>0$ in $(\ref{1.16-5})$, we arrive at
\begin{align}
B_1(t)&\geqslant-\frac{1}{2} \varepsilon^{\prime}(t)\|v\|^2-\alpha \varepsilon(t)\|u_t\|^2+\left(\frac{\alpha^2}{2} \varepsilon^{\prime}(t)-\frac{\lambda_1}{2}-\frac{\alpha}{2}\right)\|u\|^2+(\alpha-\alpha^2)\|\nabla u\|^2\non\\
&+\alpha\|\nabla u\|^{m+2}+\left\|\nabla u_t\right\|^2+\frac{\delta_1}{4}\|\eta\|_{\mathcal{M}_2}^2-C_3|\Omega|-\frac{\alpha}{2}\|h\|^2.\label{3.23-5}
\end{align}

Moreover, integrating $(\ref{3.15-5})$ from $0$ to $t$, it follows that
\begin{equation}
A_1(t)+\int_0^t B_1(s) d s \leqslant A_1(0).
\label{3.24-5}
\end{equation}

Inserting $(\ref{3.19-5})$ and $(\ref{3.23-5})$ into $(\ref{3.24-5})$, we derive
\begin{align}
&\frac{1}{2} \varepsilon(t)\|v\|^2+\frac{1}{8}\|\nabla u\|^2+\frac{1}{m+2}\|\nabla u\|^{m+2}+\frac{1}{2}\|\eta\|_{\mathcal{M}_2}^2-\frac{1}{2} \int_0^t \varepsilon^{\prime}(s)\|v(s)\|^2 d s\non \\
&-\alpha \int_0^t \varepsilon(s)\left\|u_s\right\|^2 d s+\int_0^t\left(\frac{\alpha^2}{2} \varepsilon^{\prime}(s)-\frac{\lambda_1}{2}-\frac{\alpha}{2}\right) \|u(s)\|^2 d s +\left(\alpha-\alpha^2\right) \int_0^t\|\nabla u(s)\|^2 d s\non\\
&+\alpha \int_0^t\|\nabla u(s)\|^{m+2} d s+\int_0^t\left\|\nabla u_s\right\|^2 d s+\frac{\delta_1}{4} \int_0^t\|\eta\|_{\mathcal{M}_2}^2 d s -C_3|\Omega| t-\frac{\alpha}{2}\|h\|^2 t \non\\
&\leqslant \frac{1}{2} \varepsilon(0)\|v_0\|^2+\frac{1}{8}\left\|\nabla u_0\right\|^2+\frac{1}{m+2}\left\|\nabla u_0\right\|^{m+2}+\frac{1}{2}\left\|\eta_0\right\|_{\mathcal{M}_2}^2.\label{3.25-5}
\end{align}

Hence, it follows from the Poincar\'{e} inequality that if $\alpha$ further satisfies $\left(1+\frac{L}{2 \lambda_1}\right) \alpha^2+\left(\frac{1}{2 \lambda_1}-1\right)+\frac{1}{2} \leqslant 0$ and $0<\alpha<\frac{\lambda_1}{L}$, then
\begin{align}
&-\alpha \int_0^t \varepsilon(s)\left\|u_s\right\|^2 d s+\int_0^t\left(\frac{\alpha^2}{2} \varepsilon^{\prime}(s)-\frac{\lambda_1}{2}-\frac{\alpha}{2}\right)\|u(s)\|^2 d s+\left(\alpha-\alpha^2\right)\int_0^t\|\nabla u(s)\|^2 ds\non\\
&+\int_0^t\left\|\nabla u_s\right\|^2 d s\geqslant 0.\label{3.26-5}
\end{align}

Noting that $\varepsilon(t)$ is decreasing and
\begin{equation}
\frac{1}{m+2}\|\nabla u\|^{m+2}-\frac{1}{2} \int_0^t \varepsilon^{\prime}(s)\|v(s)\|^2 d s+\alpha \int_0^t\|\nabla u(s)\|^{m+2} d s+\frac{\delta_1}{4} \int_0^t\|\eta\|^2_{\mathcal M_2}d s>0,
\label{3.27-5}
\end{equation}
then by $(\ref{3.25-5})-(\ref{3.26-5})$, we obtain
\begin{align}
\frac{1}{2} \varepsilon(t)\|v\|^2+\frac{1}{8}\|\nabla u\|^2+\frac{1}{2}\| \eta\|^2_{\mathcal M_2}&\leqslant \frac{1}{2} \varepsilon(0)\| v_0\|^2+\frac{1}{8}\| \nabla u_0\|^2+\frac{1}{m+2}\| \nabla u_0\|^{m+2} \non\\
&+\frac{1}{2}\|\eta_0\|^2_{\mathcal M_2}+C_3|\Omega| t+\frac{\alpha}{2} \| h\|^2t.\label{3.28-5}
\end{align}

Then from $v=u_t+\alpha u$, $(\ref{1.3-5})$, $(\ref{1.20-5})$, the Young and the Poincar\'{e} inequalities, we derive there exists a bounded function $R_1(t_1)$ satisfies
\begin{equation}
R_1(t_1) \geqslant C\left\{\frac{1}{2} \varepsilon(0)\|v_0\|^2+\frac{1}{8}\left\|\nabla u_0\right\|^2+\frac{1}{m+2}\|\nabla u_0\|^{m+2}+\frac{1}{2}\| \eta_0\|_{\mathcal M_2}^2+C_3 |\Omega|t+\frac{\alpha}{2}\|h\|^2 t\right\}
\label{3.29-5}
\end{equation}
for some fixed $C>0$ and a constant $C_6>\max \left\{2, \frac{2\alpha^2 L}{\lambda_1}+1\right\}$ such that
\begin{equation}
\varepsilon(t)\|u_t\|^2+\|\nabla u\|^2+\|\eta\|_{\mathcal M_2}^2 {\leq} C_6\left(\varepsilon(t)\|v\|^2+\|\nabla u\|^2+\|\eta\|_{\mathcal M_2}^2\right) \leqslant R(t_1),
\label{3.30-5}
\end{equation}
where $t_1 \ge t$.

Consequently,
\begin{equation}
\mathscr{B}_1=\left\{\left(u_0, (u_0)_t, \eta_0\right) \in \mathcal{H}\,\big|\left\|\nabla u_0\right\|^2+\varepsilon(0)\| (u_0)_t\|^2+\left\|\eta_0\right\|_{\mathcal{M}_2}^2 \leqslant R\left(t_1\right)\right\}
\label{3.31-5}
\end{equation}
is a bounded absorbing set for the semigroup $\{S(t)\}_{t \geqslant 0}$ of the weak solutions to problem $(\ref{1.10-5})-(\ref{1.12-5})$ in the weighted time-dependent space $\mathcal H$, then it is clearly that $\bar{\mathscr B}=\bigcup\limits_{t \geqslant 0} S(t) \mathscr B_1$ is also a bounded absorbing set. $\hfill$$\Box$

Next, we will prove the the following lemma, which is crucial to prove the asymptotic compactness of the semigroup associated to problem $(\ref{1.10-5})-(\ref{1.12-5})$.
\begin{Lemma}\label{lem3.2-5} Under the assumptions of $M(\cdot)$, $\varepsilon(t)$, $k(\cdot)$, $\mu(\cdot)$, $\delta$, $f$, $g$ and $h$ in $\S 1$, if
$\left(u_1, (u_1)_t, \eta_1\right)$ and $\left(u_2, (u_2)_t, \eta_2\right)$ are two weak solutions to problem $(\ref{1.10-5})-(\ref{1.12-5})$, $\bar{u}(t)=u_1(t)-u_2(t)$, $\bar{\eta}=\bar{\eta}^t=\eta_1^t-\eta_2^t$,
\begin{equation}
\widetilde{A}_1(t)=\frac{1}{2}\|\nabla \bar{u}\|^2+\frac{1}{2} \varepsilon(t)\|\bar{u}_t\|^2+\frac{1}{2}\|\bar{\eta}\|_{\mathcal M_2}^2
\label{3.32-5}
\end{equation}
and
\begin{equation}
E(t)=\|\nabla \bar{u}\|^2+\varepsilon(t)\|\bar{u}_t\|^2+\|\bar{\eta}\|^2_{\mathcal M_2},
\label{3.33-5}
\end{equation}
then there exists some positive constants $C_{\delta, L}$, $\widetilde{C}_\delta$, $\widetilde{C}_{2, \delta}$, $C_{\beta_1}$, $C_{p_1}>C_1>0$ and $C_{\widetilde A_1}>0$ depending on the parameters of their lower angle indices, and a contractive function $\Phi_T$ such that
\begin{equation}
E(t) \leqslant 2 \widetilde A_1(T) \leqslant \frac{2C_{\widetilde A_1}}{T}+\frac{2}{T} \Phi_T\left(\left(u_1, (u_1)_t, \eta_1\right),(u_2, (u_2)_t, \eta_2\right)).
\label{3.34-5}
\end{equation}
\end{Lemma}
$\mathbf{Proof.}$ From $\left(u_1, (u_1)_t, \eta_1\right)$ and $\left(u_2, (u_2)_t, \eta_2\right)$ are two weak solutions to problem $(\ref{1.10-5})-(\ref{1.12-5})$, $\bar{u}(t)=u_1(t)-u_2(t)$ and $\bar{\eta}=\bar{\eta}^t=\eta_1^t-\eta_2^t$, we obtain
\begin{align}
&\varepsilon(t) \bar{u}_{tt}-\|\nabla u_1\|^m \Delta u_1+\| \nabla u_2 \|^m \Delta u_2-\Delta \bar{u}-\Delta \bar{u}_t-\int_0^{\infty} \mu(s) \Delta \bar{\eta}(x, s) d s+\delta f((u_1)_t)\non\\
&-\delta f((u_2)_t) =g\left(u_1\right)-g\left(u_2\right)\label{3.35-5}
\end{align}
and
\begin{equation}
\bar{\eta}_t=-\bar{\eta}_s+\bar{u}_t.
\label{3.36-5}
\end{equation}

Taking $L^2(\Omega)$ inner product with $\bar u$ in $(\ref{3.35-5})$, then integrating the obtained formula from $0$ to $T$, we conclude
\begin{align}
&\int_0^T \int_{\Omega} \varepsilon(t) \bar{u}_{tt} \cdot \bar{u} d x d r-\int_0^T \int_{\Omega}\left\|\nabla u_1\right\|^m\left(\Delta u_1\right) \bar{u} d x d r+\int_0^T \int_{\Omega}\left\|\nabla u_2\right\|^m\left(\Delta u_2\right) \bar{u} d x d r \non\\
&-\int_0^T \int_{\Omega} \Delta \bar{u} \cdot \bar{u} d x d r-\int_0^T \int_{\Omega} \Delta \bar{u}_t \cdot \bar{u} d x d r-\int_0^T \int_{\Omega} \int_0^{\infty} \mu(s) \Delta \bar{\eta}(x, s) d s \cdot \bar{u} d x d r\non\\
&+\int_0^T \int_{\Omega}\left(\delta f\left(\left(u_1\right)_t\right)-\delta f\left(\left(u_2\right)_t\right)\right) \bar{u} d x d r=\int_0^T \int_{\Omega}\left(g\left(u_1\right)-g\left(u_2\right))\bar{u} d x d r.\right.\label{3.37-5}
\end{align}

Then by $(\ref{3.37-5})$ and some simple calculations, we derive
\begin{align}
&\int_0^T\|\nabla \bar{u}\|^2 d r=-\int_{\Omega} \varepsilon(T) \bar{u}(T) \bar{u}_t(T) d x+\int_{\Omega} \varepsilon(0) \bar{u}(0) \bar{u}_t(0) d x+\int_{\Omega} \int_0^T \varepsilon^{\prime}(t) \bar{u} \bar{u}_t d r d x \non\\
&+\int_0^T \varepsilon(t)\|\bar{u}_t\|^2 d r-\int_0^T\left\|\nabla u_1\right\|^m\|\nabla \bar{u}\|^2 d r-\int_0^T \int_{\Omega}\left(\left\|\nabla u_1\right\|^m-\left\|\nabla u_2\right\|^m\right) \nabla u_2 \nabla \bar{u} d x d r \non\\
&-\frac{1}{2}\|\nabla \bar{u}(T)\|^2+\frac{1}{2}\|\nabla \bar{u}(0)\|^2-\int_0^T(\bar{\eta}, \bar{u})_{\mathcal{M}_2} d r-\int_0^T \int_{\Omega}\left(\delta f\left(\left(u_1\right)_t\right)-\delta f\left(\left(u_2\right)_t\right)\right) \bar{u} d x d r \non\\
&+\int_0^T \int_{\Omega}\left(g\left(u_1\right)-g\left(u_2\right)\right) \bar{u} d x d r.\label{3.38-5}
\end{align}

Thanks to $(\ref{1.15-5})$ and the Young inequality, we arrive at there exists a constant $\alpha_{0}>0$ such that
\begin{equation}
-\int_0^T(\bar{\eta}, \bar{u})_{\mathcal M_2} d r \leqslant-\alpha_0 \int_0^T\|\bar{\eta}\|_{\mathcal M_2}^2 d r-\delta_1 \int_0^T\|\nabla\bar{u}\|^2 d r.
\label{3.39-5}
\end{equation}

Besides, multiplying $(\ref{3.35-5})$ by $ \bar{u}_t$ and integrating the obtained formula in $\Omega$, we obtain
\begin{align}
&\frac{d}{d t}\left(\frac{1}{2} \varepsilon (t)\|\bar{u}_t\|^2+\frac{1}{2}\|\nabla \bar{u}\|^2+\frac{1}{2}\|\bar{\eta}\|_{\mathcal M_2}^2\right)-\frac{1}{2} \varepsilon^{\prime}(t)\|\bar{u}_t\|^2+\int_{\Omega}\left\|\nabla u_1\right\|^m \frac{1}{2} \frac{d}{d t}(\nabla \bar{u})^2 d x \non\\
&+\int_{\Omega}\left(\left\|\nabla u_1\right\|^m-\left\|\nabla u_2\right\|^m\right) \nabla u_2 \nabla \bar{u}_t d x+\left\|\nabla\bar{u}_t\right\|^2+\left(\bar{\eta}, \bar{\eta}_s\right)_{\mathcal M_2}-\int_{\Omega}\left(g\left(u_1\right)-g\left(u_2\right)\right) \bar{u}_t d x \non\\
&+\int_{\Omega}\left(\delta f\left((u_1\right)_t)-\delta f\left((u_2\right)_t\right)) \bar{u}_t d x=0.\label{3.40-5}
\end{align}

Through the calculations similar to $(\ref{3.2-5})$, we conclude
\begin{equation}
\left(\bar{\eta}, \bar{\eta}_s\right)_{\mathcal M_2} \geqslant \frac{\delta_2}{2}\|\bar{\eta}\|_{\mathcal M_2}^2.
\label{3.41-5}
\end{equation}

Furthermore, inserting $(\ref{3.41-5})$ into $(\ref{3.40-5})$, we deduce from $(\ref{3.32-5})$ that
\begin{align}
&\frac{d}{d t} \widetilde{A}_1(t)-\frac{1}{2} \varepsilon^{\prime}(t)\|\bar{u}_t\|^2+\int_{\Omega}\left\|\nabla u_1\right\|^m \frac{1}{2} \frac{d}{dt}(\nabla \bar{u})^2 d x+\int_{\Omega}\left(\left\|\nabla u_1\right\|^m-\left\|\nabla u_2\right\|^m\right) \nabla u_2 \nabla \bar{u}_t d x \non\\
&+\|\nabla \bar{u}_t\|^2+\frac{\delta_2}{2}\| \bar{\eta}\|_{\mathcal M_2}^2+\int_{\Omega}(\delta f((u_1)_t)-\delta f((u_2)_t)) \bar{u}_t d x-\int_{\Omega}(g(u_1)-g(u_2)) \bar{u}_t d x=0.\label{3.42-5}
\end{align}

Integrating $(\ref{3.42-5})$ from $s$ to $T$, we derive
\begin{align}
&\widetilde{A}_1(T)-\frac{1}{2} \int_s^T \varepsilon^{\prime}(t)\|\bar{u}_t\|^2 d \tau+\frac{1}{2} \int_{\Omega}\left\|\nabla u_1(T)\right\|^m(\nabla \bar{u}(T))^2 d x-\frac{1}{2} \int_{\Omega}\left\|\nabla  u_1(s)\right\|^m(\nabla \bar{u}(s))^2 d x \non\\
&-\frac{m}{2} \int_{\Omega} \int_s^T\left\|\nabla u_1\right\|^{m-1} \nabla\left(u_1\right)_t(\nabla \bar{u})^2 d \tau d x+\int_s^T \int_{\Omega}\left(\left\|\nabla u_1\right\|^m-\left\|\nabla u_2\right\|^m\right) \nabla u_2 \nabla \bar{u}_t d x d \tau \non\\
&+\int_s^T\left\|\nabla \bar{u}_t\right\|^2 d \tau+\frac{\delta_2}{2} \int_s^T\|\bar{\eta}\|_{\mathcal{M}_2}^2 d \tau+\int_s^T \int_{\Omega}\left(\delta f\left(\left(u_1\right)_t\right)-\delta f\left(\left(u_2\right)_t\right)\right) \bar{u}_t d x d \tau \non\\
&-\int_s^T \int_{\Omega}\left(g\left(u_1\right)-g\left(u_2\right)\right) \bar{u}_t d x d \tau \leqslant \widetilde{A}_1(s).\label{3.43-5}
\end{align}

Then integrating $(\ref{3.43-5})$ from $0$ to $T$, we obtain
\begin{align}
&T \widetilde{A}_1(T)-\frac{1}{2} \int_0^T \int_s^T \varepsilon^{\prime}(t)\|\bar{u}_t\|^2 d \tau d r+\frac{1}{2} \int_0^T \int_{\Omega}\left\|\nabla u_1(T)\right\|^m(\nabla \bar{u}(T))^2 d x d r\non\\
&-\frac{1}{2} \int_0^T \int_{\Omega}\|\nabla u_1(s)\|^m(\nabla \bar{u}(s))^2 dxdr-\frac{m}{2} \int_0^T \int_{\Omega} \int_s^T\left\|\nabla u_1\right\|^{m-1} \nabla\left(u_1\right)_t(\nabla \bar{u})^2 d \tau d x d r\non\\
&+\int_0^T \int_s^T \int_{\Omega}\left(\left\|\nabla u_1\right\|^m-\left\|\nabla u_2\right\|^m\right) \nabla u_2 \nabla \bar{u}_t d x d \tau d r +\int_0^T \int_s^T\|\nabla \bar{u}_t\|^2 d \tau d r\non\\
&+\frac{\delta_2}{2} \int_0^T \int_s^T\|\bar{\eta}\|^2 _{\mathcal M_2}d \tau d r+\int_0^T \int_s^T \int_{\Omega}(\delta f((u_1)_t)-\delta f((u_2)_t)) \bar{u}_t d x d \tau d r\non\\
&-\int_0^T \int_s^T \int_{\Omega}(g(u_1)-g(u_2)) \bar{u}_t d x d \tau d r \leqslant \int_0^T \widetilde{A}_1(s) d r.\label{3.44-5}
\end{align}

It is obvious that
\begin{equation}
\frac{1}{2} \int_0^T \int_{\Omega}\|\nabla u_1(T)\|^m(\nabla \bar{u}(T))^2 d x d r+\int_0^T \int_s^T\|\nabla \bar{u}_t\|^2 d \tau d r+\frac{\delta_2}{2} \int_0^{T} \int_s^T\|\bar{\eta}\|_{\mathcal M_2}^2d \tau d r \geqslant 0.
\label{3.45-5}
\end{equation}

Using $(\ref{1.4-5})-(\ref{1.6-5})$, $\bar u=(u_1)_t-(u_2)_t$ and $\varepsilon(t)$ is decreasing, then by some simple calculations, we derive that there exists a constant $C_{\delta, L}>0$ depending on $\delta$ and $L$ such that
\begin{align}
&-\frac{1}{2} \int_0^T \int_s^T \varepsilon^{\prime}(t)\|\bar{u}_t\|^2 d \tau d r+\int_0^T \int_s^T \int_{\Omega}\left(\delta f\left(\left(u_1\right)_t\right)-\delta f\left(\left(u_2\right)_t\right)\right) \bar{u}_t d x d \tau d r\non\\
&\geqslant C_{\delta, L} \int_0^T \int_{s}^T \int_{\Omega} \bar{u}_t^2 d x d \tau d r\geqslant 0.\label{3.46-5}
\end{align}

Thanks to $(\ref{3.44-5})-(\ref{3.46-5})$, we deduce
\begin{align}
T \widetilde{A}_1(T)& \leqslant \int_0^T \widetilde{A}_1(s) d r+\frac{1}{2} \int_0^T \int_{\Omega}\left\|\nabla u_1(s)\right\|^m(\nabla \bar{u}(s))^2 d x d r\non\\
&+\frac{m}{2} \int_0^T \int_{\Omega} \int_s^T\|\nabla u_1\|^{m-1} \nabla\left(u_1\right)_t(\nabla \bar{u})^2 d \tau d x d r\non \\
&-\int_0^T \int_s^T \int_{\Omega}\left(\left\|\nabla u_1\right\|^m-\left|\left|\nabla u_2\right|\right|^m\right) \nabla u_2 \nabla \bar{u}_t d x d \tau d r\non\\
&+\int_0^T \int_s^T \int_{\Omega}\left(g\left(u_1\right)-g\left(u_2\right)\right) \bar{u}_t d x d \tau d r.\label{3.47-5}
\end{align}

Similarly, by $(\ref{1.4-5})-(\ref{1.6-5})$, $0<\delta<1$ and $\bar u=(u_1)_t-(u_2)_t$, we conclude there exists a constant $\widetilde{C}_\delta>0$ depending on $\delta$ such that
\begin{equation}
\int_0^T \int_{\Omega}(\delta f ((u_1)_ t)-\delta f((u_2)_t)) \bar{u}_t d x d r \geqslant \widetilde{C}_\delta+\int_0^T\|\bar{u}_t\|^2 d r.
\label{3.48-5}
\end{equation}

Moreover, integrating $(\ref{3.42-5})$ from $0$ to $T$, we arrive at
\begin{align}
&\widetilde{A}_1(T)-\frac{1}{2} \int_0^T \varepsilon^{\prime}(t)\|\bar{u}_t\|^2 d r+\int_0^T \int_{\Omega}\left\|\nabla u_1\right\|^m \frac{1}{2} \frac{d}{d t}(\nabla \bar{u})^2 d x d r\non\\
&+\int_0^T \int_{\Omega}\left(\left\|\nabla u_1\right\|^m-\left\|\nabla u_2\right\|^m\right) \nabla u_2 \nabla \bar{u}_td x d r +\int_0^T\|\nabla \bar{u}_{t}\|^2 d r\non\\
&+\frac{\delta_2}{2} \int_0^T\| \bar{\eta} \|_{\mathcal M_2}^2 d r+\int_0^T \int_{\Omega}\left(\delta f(\left(u_1\right)_t\right)-\delta f\left( u_2)_t)\right) \bar{u}_t d x d r \non\\
&-\int_0^T \int_{\Omega}\left(g\left(u_1\right)-g\left(u_2\right)\right) \bar{u}_t d x d r \leqslant \widetilde{A}_1(0).\label{3.49-5}
\end{align}

It follows from $(\ref{3.49-5})$ that
\begin{align}
&\widetilde{A}_1(T)-\frac{1}{2} \int_0^T \varepsilon^{\prime}(t)\|\bar{u}_t\|^2 d r+\frac{1}{2} \int_{\Omega}\left\|\nabla u_1(T)\right\|^m(\nabla \bar{u}(T))^2 d x \non\\
&-\frac{1}{2} \int_{\Omega}\left\|\nabla u_1(0)\right\|^m(\nabla \bar{u}(0))^2 d x-\frac{m}{2} \int_0^T \int_{\Omega}\|\nabla u_1\|^{m-1} \nabla\left(u_1\right)_t(\nabla \bar{u})^2 d x d r\non\\
&+\int_0^T \int_{\Omega}\left(\left\|\nabla u_1\right\|^m-\left\|\nabla u_2\right\|^m\right) \nabla u_2 \nabla \bar{u}_t d x d r+\int_0^T\left\|\nabla \bar{u}_t\right\|^2 d r \non\\
&+\frac{\delta_2}{2} \int_0^T\|\bar{\eta} \|_{\mathcal M_2}^2 d r+\int_0^T \int_{\Omega}(\delta f((u_1)_t))-\delta f((u_2)_t)) \bar{u}_t d x d r\non\\
&-\int_0^T \int_{\Omega}(g(u_1)-g(u_2)) \bar{u}_t d x d r \leqslant \widetilde{A}_1(0).\label{3.50-5}
\end{align}

Inserting $(\ref{3.48-5})$ into $(\ref{3.50-5})$, we derive
\begin{align}
&\widetilde{A}_1(T)-\frac{1}{2} \int_0^T \varepsilon^{\prime}(t)\|\bar{u}_t\|^2 d r+\frac{1}{2} \int_{\Omega}\left\|\nabla u_1(T)\right\|^m(\nabla \bar{u}(T))^2 d x \non\\
&-\frac{1}{2} \int_{\Omega}\left\|\nabla u_1(0)\right\|^m(\nabla \bar{u}(0))^2 d x-\frac{m}{2} \int_0^T \int_{\Omega}\|\nabla u_1\|^{m-1} \nabla\left(u_1\right)_t(\nabla \bar{u})^2 d x d r\non\\
&+\int_0^T \int_{\Omega}\left(\left\|\nabla u_1\right\|^m-\left\|\nabla u_2\right\|^m\right) \nabla u_2 \nabla \bar{u}_t d x d r+\int_0^T\left\|\nabla \bar{u}_t\right\|^2 d r+\frac{\delta_2}{2} \int_0^T\|\bar{\eta} \|_{\mathcal M_2}^2 d r\non\\
&+\widetilde{C}_\delta+\int_0^T\|\bar{u}_t\|^2 d r-\int_0^T \int_{\Omega}(g(u_1)-g(u_2)) \bar{u}_t d x d r \leqslant \widetilde{A}_1(0).\label{3.51-5}
\end{align}

Thanks to $(\ref{1.2-5})$, $(\ref{1.3-5})$, $(\ref{3.51-5})$ and $\varepsilon(t)$ is decreasing, we obtain
\begin{align}
&\widetilde{A}_1(T)+\left(1+\frac{L}{2}\right) \int_0^T\|\bar{u}_t\|^2 d r+\int_0^T\left\|\nabla \bar{u}_t\right\|^2 d r+\frac{\delta_2}{2} \int_0^T\|\bar{\eta}\|_{\mathcal M_2}^2 d r \non\\
&\leqslant-\frac{1}{2} \int_{\Omega}\|\nabla u_1(T)\|^m(\nabla \bar{u}(T))^2 d x+\frac{1}{2} \int_{\Omega}\| \nabla u_1(0)\|^m(\nabla \bar{u}(0))^2 d x-\widetilde{C}_\delta+\widetilde{A}_1(0)\non\\
&+\frac{m}{2} \int_0^T \int_{\Omega}\| \nabla u_1 \|^{m-1} \nabla\left(u_1\right)_t(\nabla \bar{u})^2 d x d r+\int_0^T \int_{\Omega}\left(g\left(u_1\right)-g(u_2\right)) \bar{u}_t d x d r\non\\
&+\int_0^{T} \int_{\Omega}\left(\left\|\nabla u_1\right\|^m-\| \nabla u_2 \|^m\right) \nabla u_2 \nabla \bar{u}_t d x d r.\label{3.52-5}
\end{align}

Then from $\widetilde{A}_1(T)+\int_0^T\left\|\nabla \bar{u}_t\right\|^2 d r \geqslant 0$ and $(\ref{3.52-5})$, we deduce there exists a constant $\widetilde{C}_{2, \delta}>0$ depending on $\delta$ such that
\begin{align}
&(1+2 L) \int_0^T\left\|\bar{u}_t\right\|^2 d r+\int_0^T\|\bar{\eta}\|_{\mathcal{M}_2}^2 d r \leqslant-\frac{\widetilde{C}_{2, \delta}}{2} \int_{\Omega}\left\|\nabla u_1(T)\right\|^m(\nabla \bar{u}(T))^2 d x\non\\
&+\frac{\widetilde{C}_{2, \delta}}{2} \int_{\Omega}\left\|\nabla u_1(0)\right\|^m(\nabla \bar{u}(0))^2 d x+ \frac{m \widetilde{C}_{2, \delta}}{2} \int_0^T \int_{\Omega}\left\|\nabla u_1\right\|^{m-1} \nabla\left(u_1\right)_t(\nabla \bar{u})^2 d xdr\non\\
&-\widetilde{C}_{2, \delta} \int_0^T \int_{\Omega}\left(\left\|\nabla u_1\left\|\left.\right|^m-\right\| \nabla u_2\right\|^m\right) \nabla u_2 \nabla \bar{u}_t d x d r+\widetilde{C}_{2, \delta} \int_0^{T} \int_{\Omega}\left(g\left(u_1\right)-g\left(u_2\right)\right) \bar{u}_t d x d r\non\\
&-\widetilde{C}_{2, \delta} \widetilde{C}_\delta+\widetilde{C}_{2, \delta}\widetilde{A}_1(0).\label{3.53-5}
\end{align}

Then inserting $(\ref{3.39-5})$ into $(\ref{3.38-5})$, we arrive at
\begin{align}
&\int_0^T\|\nabla \bar{u}\|^2 d r \leqslant-\int_{\Omega} \varepsilon(T) \bar{u}(T) \bar{u}_t(T) d x+\int_{\Omega} \varepsilon(0) \bar{u}(0) \bar{u}_t(0) d x+\int_{\Omega} \int_0^T \varepsilon^{\prime}(t) \bar{u} \bar{u}_t d r d x \non \\
&+\int_0^T \varepsilon(t)\|\bar{u}_t\|^2 d r-\int_0^T\left\|\nabla u_1\right\|^m\|\nabla \bar{u}\|^2 d r-\int_0^T \int_{\Omega}\left(\left\|\nabla u_1\right\|^m-\left\|\nabla u_2\right\|^m\right) \nabla u_2 \nabla \bar{u} d x d r \non \\
&-\frac{1}{2}\|\nabla \bar{u}(T)\|^2+\frac{1}{2}\|\nabla \bar{u}(0)\|^2-\int_0^T \int_{\Omega}\left(\delta f\left(\left(u_1\right)_t\right)-\delta f\left(\left(u_2\right)_t\right)\right) \bar{u} d x d r-\alpha_0 \int_0^T\|\bar{\eta}\|^2_{\mathcal M_2} d r\non \\
&-\delta_1 \int_0^T\|\nabla \bar{u}\|^2 d r+\int_0^T \int_{\Omega}\left(g\left(u_1\right)-g\left(u_2\right)\right) \bar{u} d x d r.\label{3.54-5}
\end{align}

By $(\ref{3.53-5})$ and $(\ref{3.54-5})$, we conclude
\begin{align}
&\left(1+\delta_1\right)\int_0^T\|\nabla \bar{u}\|^2 d r+(1+L) \int_0^T\|\bar{u}_t\|^2 d r+\left(1+\alpha_0\right) \int_0^T\| \bar{\eta}\|^2_{\mathcal M_2} d r \non\\
&\leqslant-\frac{\widetilde{C}_{2, \delta}}{2} \int_{\Omega}\left\|\nabla u_1(T)\right\|^m(\nabla \bar{u}(T))^2 d x +\frac{\widetilde{C}_{2, \delta}}{2} \int_{\Omega}\left\|\nabla u_1(0)\right\|^m(\nabla \bar{u}(0))^2 d x-\widetilde{C}_{2, \delta} \widetilde{C}_\delta\non\\
&+ \frac{m \widetilde{C}_{2, \delta}}{2} \int_0^T \int_{\Omega}\left\|\nabla u_1\right\|^{m-1} \nabla\left(u_1\right)_t(\nabla \bar{u})^2 dxdr+\widetilde{C}_{2, \delta} \int_0^{T} \int_{\Omega}\left(g\left(u_1\right)-g\left(u_2\right)\right) \bar{u}_t d x d r\non\\
&-\widetilde{C}_{2, \delta} \int_0^T \int_{\Omega}\left(\|\nabla u_1\|^m-\| \nabla u_2\|^m\right) \nabla u_2 \nabla \bar{u}_t d x d r -\int_0^T\left\|\nabla u_1\right\|^m\|\nabla \bar{u}\|^2 d r\non\\
&-\int_{\Omega} \varepsilon(T) \bar{u}(T) \bar{u}_t(T) d x+\int_{\Omega} \varepsilon(0) \bar{u}(0) \bar{u}_t(0) d x+\int_{\Omega} \int_0^T \varepsilon^{\prime}(t) \bar{u} \bar{u}_t d r d x +\widetilde{C}_{2, \delta}\widetilde{A}_1(0)\non\\
&-\int_0^T \int_{\Omega}\left(\left\|\nabla u_1\right\|^m-\left\|\nabla u_2\right\|^m\right) \nabla u_2 \nabla \bar{u} d x d r-\frac{1}{2}\|\nabla \bar{u}(T)\|^2+\frac{1}{2}\|\nabla \bar{u}(0)\|^2\non\\
&-\int_0^T \int_{\Omega}\left(\delta f\left(\left(u_1\right)_t\right)-\delta f\left(\left(u_2\right)_t\right)\right) \bar{u} d x d r+\int_0^T \int_{\Omega}\left(g\left(u_1\right)-g\left(u_2\right)\right) \bar{u} d x d r.\label{3.55-5}
\end{align}

Then from $(\ref{1.2-5})$, $(\ref{1.3-5})$, $(\ref{3.32-5})$ and $(\ref{3.55-5})$, we derive
\begin{align}
&\int_0^T \widetilde{A}_1(t) d r\leqslant-\frac{\widetilde{C}_{2, \delta}}{2} \int_{\Omega}\left\|\nabla u_1(T)\right\|^m(\nabla \bar{u}(T))^2 d x+\frac{\widetilde{C}_{2, \delta}}{2} \int_{\Omega}\left\|\nabla u_1(0)\right\|^m(\nabla \bar{u}(0))^2 d x\non\\
&+ \frac{m \widetilde{C}_{2, \delta}}{2} \int_0^T \int_{\Omega}\left\|\nabla u_1\right\|^{m-1} \nabla\left(u_1\right)_t(\nabla \bar{u})^2 d xdr+\widetilde{C}_{2, \delta} \int_0^{T} \int_{\Omega}\left(g\left(u_1\right)-g\left(u_2\right)\right) \bar{u}_t d x d r\non\\
&-\widetilde{C}_{2, \delta} \int_0^T \int_{\Omega}\left(\left\|\nabla u_1\left\|\left.\right.^m-\right\| \nabla u_2\right\|^m\right) \nabla u_2 \nabla \bar{u}_t d x d r-\widetilde{C}_{2, \delta} \widetilde{C}_\delta -\int_0^T\left\|\nabla u_1\right\|^m\|\nabla \bar{u}\|^2 d r\non\\
&-\int_{\Omega} \varepsilon(T) \bar{u}(T) \bar{u}_t(T) d x+\int_{\Omega} \varepsilon(0) \bar{u}(0) \bar{u}_t(0) d x+\int_{\Omega} \int_0^T \varepsilon^{\prime}(t) \bar{u} \bar{u}_t d r d x+\widetilde{C}_{2, \delta}\widetilde{A}_1(0)\non\\
&-\int_0^T \int_{\Omega}\left(\left\|\nabla u_1\right\|^m-\left\|\nabla u_2\right\|^m\right) \nabla u_2 \nabla \bar{u} d x d r-\frac{1}{2}\|\nabla \bar{u}(T)\|^2+\frac{1}{2}\|\nabla \bar{u}(0)\|^2\non\\
&-\int_0^T \int_{\Omega}\left(\delta f\left(\left(u_1\right)_t\right)-\delta f\left(\left(u_2\right)_t\right)\right) \bar{u} d x d r+\int_0^T \int_{\Omega}\left(g\left(u_1\right)-g\left(u_2\right)\right) \bar{u} d x d r.\label{3.56-5}
\end{align}

Thanks to Lemma $\ref{lem3.1-5}$ and $0<\delta<1$, we deduce there exists a constant $C_{\beta_1}>0$ such that% the following inequalities holds
\begin{equation}
\int_0^T \int_{\Omega} \delta f(\left(u_1)_t\right)\left(u_1\right)_t d x d r \leqslant C_{\beta_1}
\label{3.57-5}
\end{equation}
and
\begin{equation}
\int_0^T \int_{\Omega} \delta f(\left(u_2)_t\right)\left(u_2\right)_t d x d r \leqslant C_{\beta_1}.
\label{3.58-5}
\end{equation}

Besides, it follows from $(\ref{1.6-5})$ that there exists constants $C_{p_1}$ and $C_1$ satisfy $C_{p_1}>C_1>0$ such that
\begin{equation}
|f(s)|^{\frac{p_1+1}{p_1}} \leqslant|f(s)|^{\frac{1}{p_1}}|f(s)|  \leqslant C_{p_1}+C_{p_1}|f(s)||s|.
\label{3.59-5}
\end{equation}

Using the H$\rm \ddot{o}$lder inequality and $(\ref{3.57-5})-(\ref{3.59-5})$, we obtain there exists a constant $\widetilde{C}_{\delta,p_1}>0$ depends on $\delta$ and $p_1$ such that
\begin{align}
& \left | \int_0^T \delta f\left(\left(u_1\right)_t\right) \bar{u} d x d r \right| \leqslant \delta \left(\int_0^T \int_{\Omega} |f(u_1)_t|^{\frac{p_1+1}{p_1}} d x d r\right)^{\frac{p_1}{p_1+1}\left.\right.} \left(\int_0^T\int_{\Omega}|\bar{u}|^{p_1+1} d x d r\right)^{\frac{1}{p_1+1}}\non\\
&\leqslant \delta\left(\int_0^T \int_{\Omega} \left(C_{p_1}+C_{p_1}\left|f\left(\left(u_1)_t\right.\right)\right||\left(u_1\right)_t \mid\right)^{\frac{p_1+1}{p_1}} d x d r\right)^{\frac{p}{p_1+1}}\left(\int_0^T \int_{\Omega} |\bar{u}|^{p_1+1} d x d r\right)^{\frac{1}{p_1+1}}\non\\
&\leqslant \widetilde{C}_{\delta, p_1} C_{p_1} C_{\beta_1}\left(\int_0^T \int_{\Omega}|\bar{u}|^{p_1+1} d x d r\right)^{\frac{1}{p_1+1}}\label{3.60-5}
\end{align}
and
\begin{equation}
\begin{aligned}
\left | \int_0^T \int_{\Omega} \delta f((u_2)_ t) \bar{u} d x d r \right| \leqslant \widetilde{C}_{\delta , p_1} C_{p_1} C_{\beta_1}\left(\int_0^T \int_{\Omega}|\bar{u}|^{p_1+1} d x d r\right)^{\frac{1}{p_1+1}}.
\end{aligned}
\label{3.61-5}
\end{equation}

Then by $(\ref{3.60-5})$ and $(\ref{3.61-5})$, we derive
\begin{equation}
\left |\int_0^T \int_{\Omega}(\delta f((u_1)_t)-\delta f((u_2)_t)) \bar{u} d x d r \right| \leqslant 2 \widetilde{C}_{\delta, p_1} C_{p_1} C_{\beta_1}\left(\int_0^T \int_{\Omega}|\bar{u}|^{p_1+1} d x d r\right)^{\frac{1}{p_1+1}}.
\label{3.62-5}
\end{equation}

Thanks to $(\ref{3.56-5})$ and $(\ref{3.62-5})$, we arrive at
\begin{align}
\int_0^T \widetilde{A}_1(t) d r&\leqslant-\frac{\widetilde{C}_{2, \delta}}{2} \int_{\Omega}\left\|\nabla u_1(T)\right\|^m(\nabla \bar{u}(T))^2 d x+\frac{\widetilde{C}_{2, \delta}}{2} \int_{\Omega}\left\|\nabla u_1(0)\right\|^m(\nabla \bar{u}(0))^2 d x\non\\
&+ \frac{m \widetilde{C}_{2, \delta}}{2} \int_0^T \int_{\Omega}\left\|\nabla u_1\right\|^{m-1} \nabla\left(u_1\right)_t(\nabla \bar{u})^2 d xdr-\widetilde{C}_{2, \delta} \widetilde{C}_\delta +\widetilde{C}_{2, \delta}\widetilde{A}_1(0)\non\\
&-\widetilde{C}_{2, \delta} \int_0^T \int_{\Omega}\left(\left\|\nabla u_1\left\|\left.\right.^m-\right\| \nabla u_2\right\|^m\right) \nabla u_2 \nabla \bar{u}_t d x d r-\int_{\Omega} \varepsilon(T) \bar{u}(T) \bar{u}_t(T) d x\non\\
&+\widetilde{C}_{2, \delta} \int_0^{T} \int_{\Omega}\left(g\left(u_1\right)-g\left(u_2\right)\right) \bar{u}_t d x d r+\int_0^T \int_{\Omega}\left(g\left(u_1\right)-g\left(u_2\right)\right) \bar{u} d x d r\non\\
&+\int_{\Omega} \varepsilon(0) \bar{u}(0) \bar{u}_t(0) d x+\int_{\Omega} \int_0^T \varepsilon^{\prime}(t) \bar{u} \bar{u}_t d r d x -\int_0^T\left\|\nabla u_1\right\|^m\|\nabla \bar{u}\|^2 d r\non\\
&-\int_0^T \int_{\Omega}\left(\left\|\nabla u_1\right\|^m-\left\|\nabla u_2\right\|^m\right) \nabla u_2 \nabla \bar{u} d x d r-\frac{1}{2}\|\nabla \bar{u}(T)\|^2+\frac{1}{2}\|\nabla \bar{u}(0)\|^2\non\\
&+2 \widetilde{C}_{\delta, p_1} C_{p_1} C_{\beta_1}\left(\int_0^T \int_{\Omega}|\bar{u}|^{p_1+1} d x d r\right)^{\frac{1}{p_1+1}}.\label{3.63-5}
\end{align}

Noting that $0<s<t<T$,  then using $(\ref{3.47-5})$ and $(\ref{3.63-5})$, we deduce
\begin{align}
T\widetilde{A}_1(T)& \leqslant-\frac{\widetilde{C}_{2, \delta}}{2} \int_{\Omega}\left\|\nabla u_1(T)\right\|^m(\nabla \bar{u}(T))^2 d x+\frac{\widetilde{C}_{2, \delta}}{2} \int_{\Omega}\left\|\nabla u_1(0)\right\|^m(\nabla \bar{u}(0))^2 d x \non\\
&+\frac{m \widetilde{C}_{2, \delta}}{2} \int_0^T \int_{\Omega}\left\|\nabla u_1\right\|^{m-1} \nabla\left(u_1\right)_t(\nabla \bar{u})^2 d x d r-\widetilde{C}_{2, \delta} \widetilde{C}_\delta +\widetilde{C}_{2, \delta}\widetilde{A}_1(0)\non\\
&-\widetilde{C}_{2, \delta} \int_0^T \int_{\Omega}\left(\left.\left\|\nabla u_1\right\|\right.^m-\left\|\nabla u_2\right\|^m\right) \nabla u_2 \nabla \bar{u}_t d x d r-\int_{\Omega} \varepsilon(T) \bar{u}(T) \bar{u}_t(T) d x\non\\
&+\widetilde{C}_{2, \delta} \int_0^T \int_{\Omega}\left(g\left(u_1\right)-g\left(u_2\right)\right) \bar{u}_t d x d r -\frac{1}{2}\|\nabla \bar{u}(T)\|^2 +\int_{\Omega} \varepsilon(0) \bar{u}(0) \bar{u}_t(0) d x\non\\
&+\int_{\Omega} \int_0^T \varepsilon^{\prime}(t) \bar{u} \bar{u}_t d r d x-\int_0^T\left\|\nabla u_1\right\|^m\|\nabla \bar{u}\|^2 d r+\int_0^T \int_{\Omega}\left(g\left(u_1\right)-g\left(u_2\right)\right) \bar{u} d x d r \non\\
&-\int_0^T \int_{\Omega}\left(\left\|\nabla u_1\right\|^m-\left\|\nabla u_2\right\|^m\right) \nabla u_2 \nabla \bar{u} d x d r+\frac{1}{2} \int_0^T \int_{\Omega}\left\|\nabla u_1(s)\right\|^m(\nabla \bar{u}(s))^2 d x d r\non\\
&+2 \widetilde{C}_{\delta, p_1} C_{p_1} C_{\beta_1}\left(\int_0^T \int_{\Omega}|\bar{u}|^{p_1+1} d x d r\right)^{\frac{1}{p_1+1}}+\int_0^T \int_s^T \int_{\Omega}\left(g\left(u_1\right)-g\left(u_2\right)\right) \bar{u}_t d x d \tau d r \non\\
&+\frac{1}{2}\|\nabla \bar{u}(0)\|^2 -\int_0^T \int_s^T \int_{\Omega}\left(\left\|\nabla u_1\right\|^m-\left.\left||\nabla u_2\right|\right|^m\right) \nabla u_2 \nabla \bar{u}_t d x d \tau d r\non\\
&+\frac{m}{2} \int_0^T \int_{\Omega} \int_s^T\|\nabla u_1\|^{m-1} \nabla\left(u_1\right)_t(\nabla \bar{u})^2 d \tau d xd r.\label{3.64-5}
\end{align}

Consequently, if we take
\begin{align}
C_{\widetilde A_1}&=-\frac{\widetilde{C}_{2, \delta}}{2} \int_{\Omega}\left\|\nabla u_1(T)\right\|^m(\nabla \bar{u}(T))^2 d x+\frac{\widetilde{C}_{2, \delta}}{2} \int_{\Omega}\left\|\nabla u_1(0)\right\|^m(\nabla \bar{u}(0))^2 d x \non\\
&-\widetilde{C}_{2, \delta} \widetilde{C}_\delta+\widetilde{C}_{2, \delta} \widetilde{A}_1(0) -\int_{\Omega} \varepsilon(T) \bar{u}(T) \bar{u}_t(T) d x +\int_{\Omega} \varepsilon(0) \bar{u}(0) \bar{u}_t(0) d x\non\\
&-\frac{1}{2}\|\nabla \bar{u}(T)\|^2+\frac{1}{2}\|\nabla \bar{u}(0)\|^2\label{3.65-5}
\end{align}
and
\begin{align}
&\Phi_T\left(\left(u_1, (u_1)_t, \eta_1\right),\left(u_2, (u_2)_t, \eta_2\right)\right)
=\frac{m \widetilde{C}_{2, \delta}}{2} \int_0^T \int_{\Omega}\left\|\nabla u_1\right\|^{m-1} \nabla\left(u_1\right)_t(\nabla \bar{u})^2 d x d r\non\\
&-\widetilde{C}_{2, \delta} \int_0^T \int_{\Omega}\left(\left\|\nabla u_1\right\|^m-\left\|\nabla u_2\right\|^m\right) \nabla u_2 \nabla \bar{u}_t d x d r +\frac{1}{2} \int_0^T \int_{\Omega}\left\|\nabla u_1(s)\right\|^m(\nabla \bar{u}(s))^2 d x d r\non\\
&+\int_{\Omega} \int_0^T \varepsilon^{\prime}(t) \bar{u} \bar{u}_t d r d x-\int_0^T\left\|\nabla u_1\right\|^m\|\nabla \bar{u}\|^2 d r -\int_0^T \int_{\Omega}\left(\left\|\nabla u_1\right\|^m-\left\|\nabla u_2\right\|^m\right) \nabla u_2 \nabla \bar{u} d x d r\non\\
&+2 \widetilde{C}_{\delta, p_1} C_{p_1} C_{\beta_1}\left(\int_0^T \int_{\Omega}|\bar{u}|^{p_1+1} d x d r\right)^{\frac{1}{p_1+1}}+\frac{m}{2} \int_0^T\int_s^T \int_{\Omega} \left\|\nabla u_1\right\|^{m-1} \nabla(u_1)_t(\nabla \bar{u})^2d x  d \tau d r\non\\
&-\int_0^T \int_s^T \int_{\Omega}\left(\left\|\nabla u_1\right\|^m-\left\|\nabla u_2\right\|^m\right) \nabla u_2 \nabla \bar{u}_t d x d \tau d r+\widetilde{C}_{2, \delta} \int_0^T \int_{\Omega}\left(g\left(u_1\right)-g\left(u_2\right)\right) \bar{u}_t d x d r\non\\
&+\int_0^T \int_{\Omega}\left(g\left(u_1\right)-g\left(u_2\right)\right) \bar{u} d x d r+\int_0^T \int_s^T \int_{\Omega}\left(g\left(u_1\right)-g\left(u_2\right)\right) \bar{u}_t d x d \tau d r,\label{3.66-5}
\end{align}
then from $(\ref{3.33-5})$ and $(\ref{3.64-5})$, $(\ref{3.34-5})$ follows readily.
$\hfill$$\Box$
%\begin{Remark} Although there are some terms in $(\ref{3.66-5})$ that could be merged, we didn't do that, because it is helpful to verify the pullback asymptotically compactness of the process $\{S(t, \tau)\}_{t \ge \tau}$ for the weak solutions to problem $(\ref{1.10-5})-(\ref{1.12-5})$ in $\mathcal H$.
%\end{Remark}
\begin{Lemma} Under the assumptions of $M(\cdot)$, $\varepsilon(t)$, $k(\cdot)$, $\mu(\cdot)$, $\delta$, $f$, $g$ and $h$ in $\S 1$, the semigroup $\{S(t)\}_{t \geqslant 0}$ for problem $(\ref{1.10-5})-(\ref{1.12-5})$ is asymptotically compact in $\mathcal H$.
\label{lem3.3-5}
\end{Lemma}
$\mathbf{Proof.}$ From Lemma $\ref{lem3.1-5}$, we derive that there exists a constant $\sigma_1 >0$ such that
\begin{equation}
\frac{2 \widetilde C_{A_1}}{T} \leqslant \sigma_1.
\label{3.67-5}
\end{equation}

Then combined with Lemma $\ref{lem2.7-5}$, we only need to verify that $\Phi_T\left(\left(u_1, \partial_t u_1, \eta_1\right),\left(u_2, \partial_t u_2, \eta_2\right)\right)$ is a contraction function, that is,
\begin{equation}
\Phi_T\left(\left(u_1, \partial_t u_1, \eta_1\right),\left(u_2, \partial_t u_2, \eta_2\right)\right) \in \mathcal{C}(\bar{\mathscr{B}}).
\label{3.68-5}
\end{equation}

Suppose $(u_{\bar n}, \partial_t u_{\bar n}, \eta_{\bar n})$ with $\bar{n}\in\mathbb{N^+}$ is a weak solution to problem $(\ref{1.10-5})-(\ref{1.12-5})$ with respect to the initial value $(u_{\bar n_0}, \partial_t u_{\bar n_0}, \eta_{\bar n_0})\in \mathcal H$, then from the proof of Lemma \ref{lem3.2-5}, we deduce
\begin{equation}
u_{\bar n} \rightarrow u \text { weakly-star in } L^2\left(0, T ; H_2\right)
\label{3.69-5}
\end{equation}
and
\begin{equation}
\partial_t u_{\bar n} \rightarrow \partial_tu \text { weakly-star in } L^2\left(0, T ; H_1\right).
\label{3.70-5}
\end{equation}

Thanks to $(\ref{1.6-5})$, we derive $H_2 \hookrightarrow \hookrightarrow L^{p_1+1}(\Omega)$, which leads to
\begin{equation}
u_{\bar n} \rightarrow u \text { in } L^2\left(0, T ; H_1\right)
\label{3.71-5}
\end{equation}
and
\begin{equation}
\partial_t u_{\bar n} \rightarrow \partial_tu \text { in } L^{p_1+1}\left(0, T ; L^{p_1+1}(\Omega)\right).
\label{3.72-5}
\end{equation}

Then from the Aubin-Lions lemma (see \cite{Lions}), we conclude there exists a subsequence $(u_j, \partial_t u_j, \eta_j)$ of $(u_{\bar n}, \partial_t u_{\bar n}, \eta_{\bar n})$ such that
\begin{equation}
\lim _{\bar{n} \rightarrow \infty} \lim _{j \rightarrow \infty} \int_0^T \int_\Omega \varepsilon^{\prime}(t)\left(u_{\bar{n}}(t)-u_j(t)\right)\left(\partial_t u_{\bar{n}}(t)-\partial_t u_j(t)\right) d x d r=0,
\label{3.73-5}
\end{equation}
\begin{equation}
\lim _{\bar{n} \rightarrow \infty} \lim _{j \rightarrow \infty} \int_0^T \int_{\Omega}\left\|\nabla u_{\bar{n}}\right\|^{m-1}\left(\nabla \partial_t u_{\bar{n}}\right)\left(\nabla u_{\bar{n}}-\nabla u_{j}\right)^2 d x d r=0,
\label{3.74-5}
\end{equation}
\begin{equation}
\lim _{\bar{n} \rightarrow \infty}\lim_{j \rightarrow \infty} \int_0^T \int_\Omega\left(\| \nabla u_{\bar{n}}\left\|^m-\right\| \nabla u_j \|^m\right) \nabla u_j\left(\nabla\left(\partial_t u_{\bar{n}}\right)-\nabla\left(\partial_tu_j\right)\right) d x d r=0,
\label{3.75-5}
\end{equation}
\begin{equation}
\lim _{\bar{n} \rightarrow \infty} \lim _{j \rightarrow \infty} \int_0^T \int_\Omega\left\|\nabla u_{\bar{n}}\right\|^m\left(\nabla u_{\bar{n}}-\nabla u_j\right)^2 d x d r=0,
\label{3.76-5}
\end{equation}
\begin{equation}
\lim _{\bar{n} \rightarrow \infty}\lim _{j \rightarrow \infty} \int_0^{T} \int_\Omega\left(\left\|\nabla u_{\bar{n}}\right\|^m-\| \nabla u_j \|^m\right) \nabla u_j\left(\nabla u_{\bar{n}}-\nabla u_j\right) d x d r.
\label{3.77-5}
\end{equation}

In addition, by $(\ref{3.71-5})$ and $(\ref{3.72-5})$, we obtain
\begin{equation}
\lim _{\bar{n} \rightarrow \infty} \lim _{j \rightarrow \infty}\left(\int_0^T \int_\Omega\left|u_{\bar{n}}-u_j\right|^{p_1+1} d x dr\right)^{\frac{1}{p_1+1}}=0.
\label{3.78-5}
\end{equation}

Thanks to $(\ref{1.7-5})$, $(\ref{1.8-5})$ and $(\ref{3.70-5})$, noting that the sequence $\left\{\left(u_{\bar{n}}, \partial_{t}u_{\bar{n}}, \eta_{\bar{n}}\right)\right\}_{{\bar{n}}=1}^{\infty}$ is bounded in $\mathcal H$, we conclude
\begin{equation}
\lim _{\bar{n} \rightarrow \infty} \lim _{j \rightarrow \infty} \int_0^T \int_{\Omega}\left(g\left(u_{\bar{n}}\right)-g\left(u_j\right)\right)\left(\partial_t u_{\bar{n}}-\partial_t u_j\right) d x d r=0,
\label{3.79-5}
\end{equation}
and
\begin{equation}
\lim _{\bar{n} \rightarrow \infty} \lim _{j \rightarrow \infty} \int_0^T \int_{\Omega}\left(g\left(u_{\bar{n}}\right)-g\left(u_j\right)\right)\left(u_{\bar{n}}- u_j\right) d x d r=0.
\label{3.80-5}
\end{equation}

For any fixed $t$, by the Lebesgue dominated convergence theorem (see \cite{Evans}), we derive
\begin{align}
& \lim _{\bar{n} \rightarrow \infty} \lim _{j \rightarrow \infty} \int_0^T \int_s^T \int_{\Omega}\left(g\left(u_{\bar{n}}\right)-g\left(u_j\right)\right)\left(\partial_t u_{\bar{n}}-\partial_t u_j\right) d x d \tau d r \non\\
& =\int_0^T\left(\lim _{\bar{n} \rightarrow \infty} \lim _{j \rightarrow \infty} \int_s^T \int_{\Omega}\left(g\left(u_{\bar{n}}\right)-g\left(u_j\right)\right)\left(\partial_t u_{\bar{n}}-\partial_t u_j\right) d x d \tau\right) d r =0\\\label{3.81-5}
\end{align}

Hence, by $(\ref{3.73-5})-(\ref{3.81-5})$, $(\ref{3.68-5})$ holds.
$\hfill$$\Box$

Using above lemmas, we can derive the following theorem on the existence of global attractor $\mathcal A$ in $\mathcal H$.

\begin{Theorem}\label{th3.2-5}
Under the assumptions of $M(\cdot)$, $\varepsilon(t)$, $k(\cdot)$, $\mu(\cdot)$, $\delta$, $f$, $g$ and $h$ in $\S 1$, then the semigroup $\{S(t)\}_{t \geqslant 0}$ for problem $(\ref{1.10-5})-(\ref{1.12-5})$ has a unique global attractor $\mathcal A$ in the weighted time-dependent space $\mathcal H$.
\end{Theorem}
$\mathbf{Proof.}$ By Definition $\ref{def2.4-5}$, Lemmas $\ref{lem3.1-5}-\ref{lem3.3-5}$, the result holds directly. $\hfill$$\Box$

\section{Regularity of global attractor $\mathcal A$}
\ \ \ \ In this section, by dividing the weak solution $u$ to problem $(\ref{1.10-5})-(\ref{1.12-5})$ into two parts, then using the energy estimation method similar to that in \cite{qy.2} and some inequalities that in \cite{Q1, Q2, Q3}, we establish the regularity of global attractor $\mathcal A$ in the weighted time-dependent space $\mathcal H^{1}$.

\begin{Theorem}
Under the assumptions of $M(\cdot)$, $\varepsilon(t)$, $k(\cdot)$, $\mu(\cdot)$, $\delta$, $f$, $g$ and $h$ in $\S 1$, then global attractor $\mathcal A$ to the semigroup $\{S(t)\}_{t \ge 0}$ for problem $(\ref{1.10-5})-(\ref{1.12-5})$ is bounded in the weighted time-dependent space $\mathcal H^{1}$.
\label{th4.1-5}
\end{Theorem}

\noindent $\mathbf{Proof.}$ Suppose $w(t)$ is a weak solution to problem $(\ref{1.10-5})-(\ref{1.12-5})$ which satisfies $w(t)=w_1(t)+w_2(t)$ with $w_1(t)$ and $w_2(t)$ satisfy the following equations
\begin{equation}
\left\{\begin{array}{l}
\varepsilon(t)\left(w_1\right)_{t t}-\left(1+\| \nabla u\left\|^m\right) \Delta w_1-\Delta\left(w_1\right)_t-\int_0^{\infty} \mu(s) \Delta \eta_a(x, s) d s=0,\right. \\
\eta_a^t=-\eta_{a, s}+\left(w_1\right)_t, \\
\left.w_1(x, t)\right|_{\partial \Omega}=0,\,\, w_1(x, 0)=u_0(x),\,\,\left(w_1\right)_t(x, 0)=u_1(x),\\
\left.\eta_a^t(x, s)\right|_{\partial \Omega}=0, \,\, \eta_a^t(x, 0)=0, \,\, \eta_a^0(x, s)=\eta_0(x, s),
\end{array}\right.
\label{4.1-5}
\end{equation}
and
\begin{equation}
\left\{\begin{array}{l}
\varepsilon(t)\left(w_2\right)_{t t}-\left(1+\|\nabla u\|^m\right) \Delta w_2-\Delta\left(w_2\right)_t-\int_0^{\infty} \mu(s) \Delta \eta_b(x, s) d s+\delta f\left(w_t\right)=g(w)+h(x), \\
\eta_b^t=-\eta_{b, s}+\left(w_2\right)_t, \\
\left.w_2(x, t)\right|_{\partial \Omega}=0, \, w_2(x, 0)=0, \,\left(w_2\right)_t(x, 0)=0, \\
\left.\eta_b^t(x, s)\right|_{\partial \Omega}=0, \,\, \eta_b^t(x, 0)=0,\,\, \eta_b^0(x, s)=0,
\end{array}\right.
\label{4.2-5}
\end{equation}
respectively, where $\eta_t=\left(\eta_a\right)_t+\left(\eta_b\right)_t$, $\left(\eta_a\right)_t=-\eta_{a, s}+\left(w_1\right)_t$ and $\left(\eta_b\right)_t=-\eta_{b, s}+\left(w_2\right)_t$.

Choosing $-\Delta\left(w_1\right)_t$ as the test function of $(\ref{4.1-5})_1$, we derive
\begin{align}
& \frac{d}{d t}\left(\left\|\Delta w_1\right\|^2+\varepsilon(t)\left\|\nabla\left(w_1\right)_t\right\|^2+\|\eta_a \|^2_{\mathcal M_3}\right)-\varepsilon^{\prime}(t)\left\|\nabla\left(w_1\right)_t\right\|^2+\|\nabla u\|^m \frac{d}{d t}\left\|\Delta w_1\right\|^2 \non\\
& +2\left\|\Delta\left(w_1\right)_t\right\|^2+2\left(\eta_a, \eta_{a, s}\right)_{\mathcal M_3}=0.\label{5.3-5}
\end{align}

From $(\ref{1.16-5})$, we obtain
\begin{equation}
\begin{aligned}
2\left(\eta_a, \eta_{a, s}\right) & =2 \int_0^{\infty} \mu(s) \int_\Omega \Delta \eta_a \Delta \eta_{a, s} d x d s =-\frac{1}{2} \int_0^{\infty} \mu^{\prime}(s)\left\|\Delta \eta_a\right\|^2 d s\geqslant \delta_2\left\|\eta_{a}\right\|^2_{\mathcal M_3}.
\end{aligned}
\label{4.4-5}
\end{equation}

By Lemma $\ref{lem3.1-5}$, we conclude there exists a constant $R_0$ satisfies $0<R_0 \leqslant\|\nabla w\|^m \leqslant R\left(t_1\right)$ such that
\begin{align}
&\frac{d}{d t}\left(\left(1+R_0\right)\left\|\Delta w_1\right\|^2+\varepsilon(t)\left\|\nabla\left(w_1\right)_t\right\|^2+\left\|\eta_{a}\right\|_{\mathcal M_3}^2\right)-\varepsilon^{\prime}(t)\left\|\nabla\left(w_1\right)_t\right\|^2+2\left\|\Delta\left(w_1\right)_t\right\|^2\non\\
&+\delta_2\|\eta_a\|_{\mathcal M_3}^2 \leqslant 0.\label{4.5-5}
\end{align}

Thanks to $\varepsilon(t)$ is decreasing, $\delta_2\ge\delta_1\ge0$ and $R_0>0$, then there exists constants $\epsilon_a>0$ and $ \epsilon_{a_2}>0$ such that
\begin{equation}
\frac{d}{d t}\left(\left\|\Delta w_1\right\|^2+\varepsilon(t)\left\|\nabla\left(w_1\right)_t\right\|^2+\|\eta_a\|_{\mathcal M_3}^2\right)+\epsilon_a\left(\left\|\Delta w_1\right\|^2+\varepsilon(t)\left\|\nabla\left(w_1\right)_t\right\|^2+\left\|\eta_a\right\|_{\mathcal M_3}^2\right) \leqslant \epsilon_{a_2}^2.
\label{4.6-5}
\end{equation}

Using the Gronwall inequality, we deduce
\begin{equation}
\left\|w_1(t)\right\|_{\mathcal H^{1}}^2 \leqslant\left\|w_1(0)\right\|_{\mathcal H^{1}}^2 e^{-\epsilon_at}+\epsilon_{a_3},
\label{4.7-5}
\end{equation}
where $\epsilon_{a_3}=\frac{\epsilon_{a_2}^2}{\epsilon_a}>0.$

Similarly, choosing $-\Delta\left(w_2\right)_t$ as the test function of $(\ref{4.2-5})_1$, we conclude
\begin{align}
&\frac{d}{d t}\left(\left\|\Delta w_2\right\|^2+\varepsilon(t)\left\|\nabla\left(w_2\right)_t\right\|^2+\left\|\eta_b\right\|_{\mathcal M_3}^2\right)-\varepsilon^{\prime}(t)\left\|\nabla\left(w_2\right)_t\right\|^2+\|\nabla u\|^m \frac{d}{d t}\left\|\Delta w_2\right\|^2 \non\\
&+2\|\Delta(w_2)_t\|^2+2(\eta_{b,} \eta_{b, s})_{\mathcal M_3}+2 \delta(f(w_t),-\Delta(w_2)_t)\non\\
&= 2(g(w),-\Delta(w_2)_t)+2(h(x),-\Delta(w_2)_t).\label{4.8-5}
\end{align}

By $(\ref{1.4-5})-(\ref{1.6-5})$, $0<\delta<1$ and the Young inequality, we obtain there exists a constant $\widetilde{C}_1>0$ such that
\begin{equation}
2| \delta\left(f\left(w_t\right),-\Delta\left(w_2\right)_t\right) \mid \leqslant \delta^2 \widetilde{C}_1\left(1+\left\|\left.w_t\|^{p_1}\right)+\right\| \Delta\left(w_2\right)_t \|^2\right..
\label{4.9-5}
\end{equation}

Besides, it follows from the Young inequality that
\begin{equation}
2\left|\left(h(x),-\Delta\left(w_2\right)_t\right)\right| \leqslant 2\|h\|^2+\frac{1}{2}\left\|\Delta\left(w_2\right)_t\right\|^2.
\label{4.10-5}
\end{equation}

By some calculations similar to $(\ref{4.4-5})$, we arrive at
\begin{equation}
2\left(\eta_b, \eta_{b, s}\right) \geqslant \delta_2\left\|\eta_b\right\|_{\mathcal M_3}^2.
\label{4.11-5}
\end{equation}

Inserting $(\ref{4.9-5})-(\ref{4.11-5})$ into $(\ref{4.8-5})$, we deduce
\begin{align}
& \frac{d}{d t}\left(\left\|\Delta w_2\right\|^2+\varepsilon(t)\left\|\nabla\left(w_2\right)_t\right\|^2+\left\|\eta_b\right\|^2_{\mathcal M_3}\right)-\varepsilon^{\prime}(t)\left\|\nabla\left(w_2\right)_t\right\|^2+\|\nabla u\|^m \frac{d}{d t}\left\|\Delta w_2\right\|^2 \non\\
& +2\left\|\Delta\left(w_2\right)_t\right\|^2+\delta_2\left\|\eta_b\right\|_{\mathcal M_3}^2+\delta^2 \widetilde{C}_1\left(1+\left\|\left.w_t \|^{p_1}\right)+\right\| \Delta\left(w_2\right)_ t \|^2\right. \non\\
& \leqslant \lambda_1^2\left\|w\right\|^2+\frac{3}{2}\left\|\Delta\left(w_2\right)_t\right\|^2+2\|h\|^2.\label{4.12-5}
\end{align}

Then from Lemma $\ref{lem3.1-5}$, we deduce that there exists constants $\epsilon_b>0$ and $C_{R\left(t_1\right)}>0$ depending on $\varepsilon(t)$, $m$, $\delta$, $\widetilde{C}_{1}$, $p_1$ $\lambda_1$ and $h$ such that
\begin{align}
& \left.\left.\frac{d}{dt}\left(\left\|\Delta w_2\right\|^2+\varepsilon(t)\left\|\nabla\left(w_2\right)_t\right\|^2+\|\eta_b\|_{\mathcal M_3}^2\right)+\epsilon_b\left(\left\|\Delta w_2\right\|^2+\varepsilon(t)\right.\left\|\nabla\left(w_2\right)_t\right\|^2+\|\eta_b\|_{\mathcal M_3}^2\right.\right) \non\\
& \leqslant C_{R\left(t_1\right)}^2. \label{4.13-5}
\end{align}

Using the Gronwall inequality, we obtain
\begin{equation}
\left\|w_2(t)\right\|_{\mathcal{H}^1}^2 \leqslant\left\|w_2(0)\right\|_{\mathcal{H}^{1}}^2 e^{-\epsilon_bt}+\frac{C_{R\left(t_1\right)}^2}{\epsilon_b}.
\label{4.14-5}
\end{equation}

By $(\ref{4.7-5})$ and $(\ref{4.14-5})$, we conclude that there exists a family
\begin{equation}
\bar{\mathscr B}_2\left(R_2\right)=\left\{w(t) \in \bar{\mathscr B}_2:\|w(t)\|_{\mathcal H^{1}}^2 \leqslant R_2\right\}
\label{4.15-5}
\end{equation}
such that
\begin{equation}
\lim _{t \rightarrow+\infty} \operatorname{dist}\left(\mathcal{A}, \bar{\mathscr B}_2\left(R_2\right)\right)=0,
\label{4.16-5}
\end{equation}
for any $0\le t \in \mathbb R$.

Therefore, we conclude that $\mathcal{A} \subseteq {\bar{\mathscr B}}(R_2)$, that is, global attractor $\mathcal A$ is bounded in the weighted time-dependent space ${\mathcal H}^1$. $\hfill$$\Box$

$\mathbf{Acknowledgment}$

This paper was supported by the China Scholarship Council with number 202206630048, the National Natural Science Foundation of China with contract number 12171082, the fundamental research funds for the central universities with contract numbers $2232022G$-$13$, $2232023G$-$13$ and a grant from science and technology commission of Shanghai municipality.

$\mathbf{Conflict\,\,of\,\,interest\,\,statement}$

The authors have no conflict of interest.

\end{document}